\documentclass[a4paper,11pt,reqno]{amsart}

\usepackage[hidelinks]{hyperref}
\usepackage[centering, top=3cm, bottom=2.5cm]{geometry}

\usepackage[foot]{amsaddr}

\usepackage{calc}
\newsavebox\CBox
\newcommand\hcancel[2][0.5pt]{%
  \ifmmode\sbox\CBox{$#2$}\else\sbox\CBox{#2}\fi%
  \makebox[0pt][l]{\usebox\CBox}%  
  \rule[0.5\ht\CBox-#1/2]{\wd\CBox}{#1}}
\usepackage{blindtext}
\usepackage{color}
\usepackage{hyperref,url}
\usepackage{float}
\usepackage{hyperref}
\usepackage{enumerate}
\usepackage{enumitem}   
\usepackage{bbm}
\usepackage{cancel}
\usepackage{mathrsfs}
\usepackage{colonequals}
\newcommand\smallO{
  \mathchoice
    {{\scriptstyle\mathcal{O}}}% \displaystyle
    {{\scriptstyle\mathcal{O}}}% \textstyle
    {{\scriptscriptstyle\mathcal{O}}}% \scriptstyle
    {\scalebox{.7}{$\scriptscriptstyle\mathcal{O}$}}%\scriptscriptstyle
  }

\usepackage{amsmath,amsfonts,amsthm,bm}
\usepackage{mathrsfs}  
\usepackage{xcolor}
\usepackage{amsfonts}
\usepackage{amssymb}
\usepackage{amsmath}
\usepackage{mathtools}
\usepackage{mathrsfs}  

\usepackage[normalem]{ulem}

\usepackage[toc,page]{appendix}
\usepackage{tikz-cd}
\theoremstyle{definition}
\newtheorem{definicao}{Definition}[section]

\theoremstyle{plain}

\theoremstyle{plain}

\newtheorem{teorema}[definicao]{Theorem}
\newtheorem{suposicao}{Hypothesis}
\newtheorem{remark}[definicao]{Remark}

\newtheorem{lema}[definicao]{Lemma}
\newtheorem{exemplo}[definicao]{Example}
\newtheorem{proposicao}[definicao]{Proposition}
\newtheorem{corolario}[definicao]{Corollary}
\newtheorem{notacao}[definicao]{Notation}

\theoremstyle{remark}

\newtheorem {mtheorem} {\bf Theorem}

\newtheorem{manualtheoreminner}{\bf Claim}

\newtheorem{manualtheoreminnerb}{\bf Step}
\newenvironment{step}[1]{%
  \renewcommand\themanualtheoreminnerb{#1}%
  \manualtheoreminnerb
}{\endmanualtheoreminner}

\definecolor{ao(english)}{rgb}{0.0, 0.5, 0.0}
\definecolor{dmagenta}{RGB}{139, 0, 139}
\definecolor{dgreen}{RGB}{0,90,0}
\definecolor{navy}{RGB}{0,0,128}

\usepackage{stackengine}

\def\d{\mathrm d}
\def \d {\mathrm{d}}
\def \supp {\mathrm{supp}}

\def \s {\mathrm{span}}
\def \cP {\mathcal{P}}
\def \cL {\mathcal{L}}
\def\fa{\text{for all}}
\definecolor{iblue}{RGB}{0, 35, 194}

\author[M.M.~Castro]{Matheus M. Castro$^{1}$}
\author[J.S.W.~Lamb]{Jeroen S.W. Lamb$^{1,2,3}$}
\author[G.~Olicón-Méndez]{Guillermo Olicón-Mendez$^{4}$}
\author[M.~Rasmussen]{Martin Rasmussen$^{1}$}

\email{\href{mailto:m.de-castro20@imperial.ac.uk}{m.de-castro20@imperial.ac.uk}}
\email{\href{mailto:g.olicon.mendez@fu-berlin.de}{g.olicon.mendez@fu-berlin.de}}
\email{\href{mailto:jeroen.lamb@imperial.ac.uk}{jeroen.lamb@imperial.ac.uk}}
\email{\href{mailto:m.rasmussen@imperial.ac.uk}{m.rasmussen@imperial.ac.uk}}

\address{$^{1}$Department of Mathematics, Imperial College London, London SW7 2AZ, UK}
\address{$^{2}$International Research Center for Neurointelligence, The University of Tokyo, Tokyo 113-0033, Japan}
\address{$^{3}$Centre for Applied Mathematics and Bioinformatics, Department of Mathematics and Natural Sciences, Gulf University for Science and Technology, Halwally 32093, Kuwait}
\address{$^{4}$ Institut für Mathematik, Freie Universität Berlin, Arnimallee 7 14195 Berlin, Germany}

\title[Quasi-stationary and quasi-ergodic measures for absorbing Markov chains]
{Existence and Uniqueness of Quasi-stationary and Quasi-ergodic Measures for Absorbing Markov chains:\\ A Banach Lattice Approach }
\allowdisplaybreaks
\begin{document}



\maketitle
% 
% \begin{align*}
%    \small {}^1 \textit{Imperial College London} \quad \quad  {}^2 \textit{The University of Tokyo}\\ {}^4 \textit{Freie Universität Berlin}
% \end{align*}

\begin{abstract}
% We study the existence and uniqueness of quasi-stationary and quasi-ergodic measures for an almost surely absorbing  Markov chain $X_n$. We find { hypotheses that depend only on the {evolution of the paths} and regularity of the transition function of $X_n$, which implies that} the two afford mentioned measures exist, and they are unique. {Our results substantially extend the conditions in which quasi-stationary and quasi-ergodic measures are known to exist. The techniques used to prove these results exploit} the Banach lattice properties of the transition function of the Markov chain $X_n.$

We establish the existence and uniqueness of quasi-stationary and quasi-ergodic measures for almost surely absorbed discrete-time Markov chains under weak conditions. We obtain our results by exploiting Banach lattice properties of transition functions under natural regularity assumptions.

\end{abstract}

%\tableofcontents
\section{Introduction and motivation}

The existence and uniqueness of quasi-stationary and quasi-ergodic measures is a central question in the theory of absorbing Markov processes, but sufficient conditions have been tailored to specific contexts such as stochastic differential equations  \cite{benaim2021,Champ,Champ2,QSB,Esperanca,Ale}. In this paper, we prove the existence and uniqueness of quasi-stationary and quasi-ergodic measures for absorbing Markov chains under weak continuity and irreducibility assumptions. Our results substantially extend the settings in which quasi-stationary and quasi-ergodic measures are known to exist, including, in particular, random systems with bounded noise.

We use Banach lattices \cite{B2,B1,Ideal1} to address the problem of the existence and uniqueness of quasi-stationary and quasi-ergodic measures from a functional-analytic point of view by considering transition probabilities as a bounded linear operator. Banach lattices allow us to study its spectrum and subsequently construct the desired measures.
Previously, Banach lattice theory has been employed to construct stationary measures and obtain statistical results for conservative Markov processes without absorption (see \cite{ CarlangeloTransf,Op, foguel,UK}). In this paper, we extend the use of this theory to absorbing Markov chains.

We consider a Markov chain $(X_{n})_{n\in\mathbb N}$ on a metric space $E$, and study the behaviour of this chain conditioned on survival in some compact set $M\subset E$. Recall that a Borel measure $\mu$ on $M$ is called \emph{quasi-stationary measure} for $X_n$ if for every $n\in \mathbb N$ and measurable set $A\subset M$, we have
$$\mathbb P_{\mu}[X_n \in A \mid \tau >n] := \frac{\int_{M} \mathbb P_x[X_n \in A] \mu(\d x)}{\int_{M} \mathbb P_x[X_n \in M] \mu(\d x)} = \mu(A)\,,$$
where $\tau$ is the stopping time for $X_n$ conditioned upon survival in $M$.

Quasi-stationary measures generalise stationary measures \cite[Section~5]{Ale} in the sense that they allow for escape from the state space. They have received increased attention in recent years \cite{Champ, Champ2,QSB}. We refer to \cite{survey} for a bibliography on the topic. Quasi-stationary measures may be used for a variety of purposes, including determining average escape times (see \cite[Lemma 5.6]{Ale}). In addition, in several contexts, it can be shown that the so-called \emph{Yaglom limit},
\begin{align}
    \lim_{n\to\infty} \mathbb P_x \left[X_n\in A \mid \tau >n\right] := \lim_{n\to\infty} \frac{\mathbb P_x \left[X_n\in A\right]}{\mathbb P_x \left[X_n\in M\right]} \label{limite},
\end{align}
converges to the unique quasi-stationary measure  \cite{Yaglom2, Yaglom1,Yaglom3}. A commonly employed technique for establishing the existence and uniqueness of quasi-stationary measures for $X_n$, as well as the existence of the Yaglom limit, relies on showing that the transition kernel of the process can be interpreted as a compact operator on a suitable space, together with the presence of a spectral gap (see \cite{benaim2021, collet2,new, QEM1}). In this paper, we provide easy-to-check sufficient conditions that guarantee these properties.

In contrast to the ergodic stationary case, the quasi-stationary measure is not the correct object to employ when studying Birkhoff averages for absorbed processes. Although quasi-stationary measures are relevant for statistics in terms of the limit  \eqref{limite}, one needs to introduce the notion of quasi-ergodic measures in order to compute conditioned Birkhoff averages. Quasi-ergodic measures are probability measures $\eta$ on $M$ that satisfy 
\begin{align}
     \lim_{n\to\infty}\mathbb E_x \left[\frac{1}{n} \sum_{i=0}^{n-1} g \circ X_i \hspace{0.1cm}\bigg\vert\hspace{0.1cm} \tau >n\right]= \int_M g(y) \eta(\d y)\ \text{for}\ \eta\text{-a.e.}~x\in M\label{piano1} 
\end{align}
for every bounded measurable function $g:M\to \mathbb R$.

While quasi-stationary measures have been well-studied, less is known about quasi-ergodic measures. Quasi-ergodic limits of the form \eqref{piano1} were first considered by Darroch and Seneta \cite{d}, who established their existence for finite irreducible Markov chains.  Their study was followed by Breyer and Roberts \cite{QEM} and later by Champagnat and Villemonais \cite{Champ,Champ2} to obtain conditions that guarantee the existence of a quasi-ergodic measure for Markov processes defined in a general state space, also providing a description of the quasi-ergodic measure in relation to the quasi-stationary measure. Explicit formulas for quasi-ergodic measures for reducible finite absorbing Markov chains have been obtained in \cite{martin1}. Quasi-ergodic measures are promising tools for analysing random dynamical systems; for instance, they are crucial for the existence of so-called conditioned Lyapunov exponents~\cite{LExp}.

Recently, numerous contributions have been made to understand the Yaglom limit \cite{R1, DifR5, DifR6, DifR8, FiniteR4}.  In particular, in \cite{Champ}, necessary and sufficient conditions for quasi-stationary measures have been established, addressing the issue of uniform exponential convergence. Although these conditions are sharp in the context of uniform exponential convergence, verifying them in applications outside the stochastic differential equations framework is complicated due to their abstract formulation. Moreover, they rely on the standing assumption that the probability of remaining in the allowed state space $M$ is positive for all positive times and all initial conditions $x\in M$. Notice that typically, this assumption is not fulfilled on \textit{repelling sets} of random maps with bounded noise (see Example~\ref{EXAMPLE:boundednoise}).

 Additionally to the uniform convergence scenario,   \cite{Champ2,kolb2,van,Kolb1}, have provided valuable insights. These studies have examined cases where the \eqref{limite} converges to a unique quasi-stationary measure in the total variation norm, although this convergence is not uniform in $x$.  This is relevant since uniform convergence of \eqref{limite} in the total variation norm is often absent in random dynamical systems with bounded noise; see Subsection~\ref{sec:RDS} below. The existence of quasi-stationary and quasi-ergodic measures in such settings is established in this paper.

As in \cite{Champ,Champ2}, alternative conditions for the existence of quasi-stationary measures found in the literature are also aimed at such stochastic differential equations \cite{Dennis1,Example10,L,V2}. These conditions are based either on lower-bound estimations of $\mathbb P_x[X_n \in \cdot \mid \tau > n]$ (such as the Doeblin condition) or on the existence of a Lyapunov function for the transition kernel. Even in elementary examples, these conditions are usually difficult to verify, and their dynamical interpretation is often unclear. In light of this, we establish the existence of both quasi-stationary and quasi-ergodic measures under natural, weak and easily verifiable dynamical conditions.

This paper is divided into six sections and one appendix.
In Section \ref{main}, we briefly recall the basic concepts of the theory of absorbing Markov chains. We introduce the main underlying hypothesis of this paper (Hypothesis \ref{(H)}), we state the main results of this paper (Theorems \ref{thm:QSM}, \ref{thm:Spec} and \ref{thm:QEM}), and apply them to some fundamental examples.
In Section \ref{sec:P}, we derive some direct consequences of Hypothesis \ref{(H)}. 
Section \ref{4} is dedicated to a brief presentation of Banach lattice theory, and the proof of Theorem~\ref{3o}. This is central to the proof of our main results.
In Section \ref{QSM1}, we show the existence and uniqueness of quasi-stationary measures for Markov chains that fulfil  Hypothesis \ref{(H)}, and prove Theorem \ref{thm:QSM}. Finally, in Section \ref{6}, we prove the existence and uniqueness of quasi-ergodic measures under Hypothesis \ref{(H)}, and we also prove Theorems \ref{thm:Spec} and \ref{thm:QEM}. \ref{appendixA} contains the proof of  Lemma~\ref{Lobo}, a technical yet essential result for the existence and uniqueness of quasi-ergodic measures.

\section{Main results \label{main}}

Let $(E,d)$ be a metric space and consider a compact subset $M\subset E$ endowed with the induced topology. We assume that $(M,\mathscr{B}(M),\rho)$ is a Borel finite measure space, where $\mathscr{B}(M)$ denotes the Borel $\sigma$-algebra of $M$. Throughout this paper, we aim to study Markov chains that are killed when the process escapes the region $M$. Since the behaviour of $X_n$ in the set $E\setminus M$ is not relevant for the desired analysis, we identify $E\setminus M$ as a single point $\partial$ and consider ${E_M} = M \sqcup \{\partial\}$ as the topological space generated by the topological basis
$$\beta = \{B_r (x);\ x\in M \ \text{and }r\in \mathbb R\}\cup \{\partial\}, $$
where $B_r(x) := \{y\in M; d(x,y) <r\},$ and $\sqcup$ denotes disjoint union.

In this paper, we assume that $$X:=\left(\Omega, \left\{\mathcal F_n\right\}_{n\in\mathbb N_0}, \left\{X_n\right\}_{n\in\mathbb N_0}, \left\{\mathcal P^n\right\}_{n\in\mathbb N_0}, \{\mathbb P_x\}_{x \in {{E_M}}}\right)$$ is a  Markov chain with state space ${E_M}$, in the sense of \cite[Definition III.1.1]{RW}. This means that the pair $(\Omega,\{\mathcal F_n\}_{n\in\mathbb N})$ is a filtered space; $X_n$ is an $\mathcal F_n$-adapted process with state space ${E_M}$; $\mathcal P^n$ a time-homegenous transition probability function of the process $X_n$ satisfying the usual measurability assumptions and Chapman-Kolmogorov equation; $\{\mathbb P_x\}_{x\in{E_M}}$ is a family of probability function satisfying $\mathbb P_x[X_0=x] = 1$ for every $x\in {E_M}$; and for all $m,n\in \mathbb N_0$, $x\in E_M$, and  every bounded measurable function $f$ on ${E_M}$, 
$$\mathbb E_x\left[f\circ X_{m+n}\mid \mathcal F_n\right] = ({\mathcal P}^m f)(X_n)\   \mathbb P_x\text{-almost surely}.$$

{As mentioned before}, we assume that $X_n$ is a Markov chain that is absorbed at $\partial,$ i.e. ${\mathcal P}(\partial,\{\partial\}) = 1. $  Given the nature of the process $X_n$ {it is natural to define the} stopping time $\tau(\omega) := \inf\{n\in\mathbb N; X_n(\omega) \not\in M\}.$
% {The above stopping time will be important during the whole paper, since the definitions of quasi-stationary and quasi-ergodic measure depend on $\tau$. In particular, observe that using  the Markov property one can achieve
% $$\mathbb P_x \left[\tau >n\right] = \mathbb E_x\left[\mathbbm 1_M \circ X_n\right]= \mathcal P^n(x,M). $$}

 We introduce some notation that is used throughout the paper. 

\begin{notacao}
Given a measure $\mu$ on $M$, we set $\mathbb P_\mu (\cdot) := \int_{M} \mathbb P_x (\cdot) \mu (\d x).$

We consider the set $\mathcal F_b(M)$ as the set of bounded Borel measurable functions on $M$.  Given $f\in {\mathcal F_b}(M)$ write
$$\mathcal P^n(f)(x):= \mathcal P^n\left(\mathbbm 1_M   f \right)(x) = \int_M f(y) \mathcal P^n(x, \mathrm{d}y), $$
\text{and} $f\circ X_n := \mathbbm 1_M \circ X_n \ f\circ X_n.$

 For every $p\in[1,\infty] $ we  denote de space $L^p \left(M,{\mathscr B}(M),\rho\right)$ by $L^p(M);$ $\mathcal C^0(M)$ refers to the set of continuous functions $f:M\to M$; and $\mathcal M(M)$ the set of Borel signed-measures on $M.$

Finally, we define the sets
$\mathcal C^0_+(M) := \{f \in \mathcal C^0(M); f\geq 0\}$  and $\mathcal M_+(M) := \{\mu\in\mathcal M(M);\ \mu(A)\geq 0,\ \text{for every }A\in\mathscr B(M)\}.$

\end{notacao}

{In the following, we recall the definition of a quasi-stationary measure.}
\begin{definicao}
A Borel measure $\mu$ on $M$ is said to be a \emph{quasi-stationary measure} for the Markov chain $X_n$ if $\mathbb P_\mu\left[X_n \in \cdot \mid \tau >n\right] =\mu(\cdot), \ {\   \text{for all}} \ n\in\mathbb N.$
We call $\lambda = \int_M\mathcal P(x,M)\mu(\d x)$ {the} \emph{survival rate} of $\mu$.
\end{definicao}

\begin{remark}
Note that since $\{\partial\}$ is absorbing
$$\mathbb P_\mu\left[X_n \in \cdot \mid \tau >n\right]  = \frac{\int_M \mathcal P^n(x,\cdot) \mu(\d x)}{\int_M \mathcal P^n(x,M) \mu(\d x)}, $$
for every $\mu \in \mathcal M_+(M).$
\end{remark}

Our initial goal is to establish the existence of quasi-stationary and quasi-ergodic measures for a  Markov chain $X_n$. We mention that our results also cover the case where $X_n$ has almost surely escaping points, i.e. points $x\in M,$ such that $\mathcal P^k(x,M) = 0$ for some $k\in \mathbb N.$ This occurs naturally in random iterated functions with bounded noise (see Section~\ref{sec:RDS}). 

The sufficient conditions to reach our results are introduced in the form of Hypothesis \ref{(H)}. Briefly, we say that a Markov chain $X_n$ satisfies Hypothesis \ref{(H)} if its transition kernels $\mathcal P(x,dy)$ are well behaved with respect to a fixed Borel measure $\rho$ on $M$, and if $X_n$ eventually reach any open set of $M$ with positive probability provided $X_0 = x$ is a non-escaping point.

\begingroup
\renewcommand\thesuposicao{(H)} 
\begin{suposicao} \label{(H)} We say that the  Markov chain $X_n$ on $E_M$ absorbed at $\partial$ fulfils \emph{Hypothesis (H)} if the following two properties hold:
\begin{enumerate}
    \item[$\boldsymbol{\mathrm{(H1)}}$] \label{(H1)} For all $x\in M$, $\mathcal P(x,\d y) \ll \rho( \d y)$; and for every $\varepsilon>0,$ there exists $\delta>0,$ such that
    $$d(x,z)<\delta \Rightarrow \|g(x,\cdot)- g(z,\cdot)\|_1 := \int_M \left|g(x,y) - g(z,y)\right| {\rho(\d y)} <\varepsilon, $$
    where $g$ is the Radon-Nikodym derivative
    $$g(x,\cdot):= \frac{\mathcal P(x,\d y)}{{\rho(\d y)}} \in L^1(M,\rho).$$

    \item[$\boldsymbol{\mathrm{(H2)}}$] \label{(H2)} Let $Z := \left\{x\in M;\ \mathcal P^k(x,M) = 0\ \text{for some }k\in\mathbb N \right\}\,.$ Then $0< \rho(M\setminus Z)$ and for any $x\in M\setminus Z$ and open set $A\subset M\setminus Z$ in the induced topology of $M\setminus Z$, there exists a natural number $n = n(x,A)$ such that
    $$\mathbb P_x\left[X_n \in A\right] = \mathcal P^n(x,A) >0. $$
\end{enumerate}
\end{suposicao}

{We now state the first main result of this paper, asserting that Hypothesis \ref{(H)} implies the existence and uniqueness of a quasi-stationary measure for $X_n$ on $M$.}

\begin{mtheorem}\label{thm:QSM} Let $X_n$ be a  Markov chain on $E_M = M\sqcup \{\partial\}$ absorbed at $\partial$ satisfying Hypothesis \ref{(H)}, then
\item[$(a)$] if $\mathcal P(x,M)=1$ for every $x\in M$, 
 then $X_n$ admits a unique stationary probability measure $\mu$ and $\supp(\mu) = M$;
    \item[$(b)$] if there exists $x \in M\setminus Z,$ such that $\mathcal P(x,M) <1,$ then $\lim_{n\to\infty} \mathcal P^n(y, M) = 0$, for all  $y\in M,$ and the {process} $X_n$ admits a unique quasi-stationary measure $\mu$ with \begin{equation}
        \supp(\mu)\supset \bigcup_{k\in \mathbb N}\bigcup_{x\in M\setminus Z}\supp(\mathcal P^k(x,\d y)) \supset M\setminus Z\label{support}
    \end{equation} 
    and survival rate $\lambda>0$.
\end{mtheorem}

Theorem \ref{thm:QSM} is proved in Section \ref{QSM1}. 

Results, analogous to Theorem \ref{thm:QSM}, were previously obtained for families of diffeomorphisms \cite{Araujo,Kifer}. Theorem \ref{thm:QSM} extend these, as we also cover families of endomorphisms.

% We mention that results analogous to part $(a)$ in Theorem \ref{thm:QSM} can be found for Markov chains generated by a family of diffeomorphisms in \cite{Araujo,Kifer}. Our results are an extension when considering the scenario where the Markov chain possesses a unique ergodic component since our results are also valid when considering a family of endomorphisms.

The main techniques presented in this paper rely on the analysis of the spectral properties of the transition function $\mathcal P,$ seen as a linear operator (see Section~\ref{sec:P}).  We summarize the main properties of $\mathcal P$ in Theorem \ref{thm:Spec}.

\begin{mtheorem}\label{thm:Spec} Let $X_n$ be a  Markov chain on $E_M$ absorbed at $\partial$ satisfying Hypothesis \ref{(H)} and $\lambda$ be the survival rate given by Theorem \ref{thm:QSM}. Then, 
\begin{itemize}
    \item[$a)$] the stochastic Koopman operator 
    \begin{align*}
    \mathcal P: (\mathcal C^0(M),\|\cdot\|_{\infty}) &\to (\mathcal C^0(M),\|\cdot\|_{\infty})\\
    f&\mapsto  \int_M f(y) \mathcal P(x,\d y)
\end{align*}
is a compact linear operator with spectral radius $r(\mathcal P)= \lambda$;
\item[$b)$] the set of eigenvalues of $\mathcal P$ with modulus $\lambda$ is given by $\{\lambda e^{\frac{2\pi i j}{m}}; j = \{0,\ldots,m-1\}\}$ for some $m\leq \#\{\text{connected components of }M\setminus Z\}$;
\item[$c)$] for every $j\in\{0,\ldots,m-1\},$ $\mathrm{\dim}(\ker(\mathcal P - \lambda e^{\frac{2\pi i j}{m}}))=1;$ and
\item[$d)$] there exists a non-negative continuous function $f$ such that $\mathcal P f = \lambda f$ and $\{x\in M; f(x)>0\}=M\setminus Z.$
\end{itemize}
\end{mtheorem}
Theorem \ref{thm:Spec} is proved in Section \ref{ProofC}. 

{
\begin{remark}
The inequality $m\leq \#\{\text{connected components of }M\setminus Z\},$ in the above theorem, shows that the spectrum of $\mathcal P$ presents topological obstructions. Moreover, it is shown in Example~\ref{example22} that it is possible for $m$ to be strictly smaller than the number of connected components of $M\setminus Z$. 
\end{remark}
We conclude this section by stating the final main result of this paper, which concerns the existence and characterisation of quasi-ergodic measures for a Markov chain $X_n$ satisfying Hypohetesis \ref{(H)}. Recall the formal definition of a quasi-ergodic measure.}
\begin{definicao}
A measure $\eta$ is called a \emph{quasi-ergodic measure} on $M$ if, for every $h\in {\mathcal F_b}(M),$
\begin{equation}
    \lim_{n\to\infty}\mathbb E_x\left[\frac{1}{n} \sum_{i=0}^{n-1} h\circ X_i \hspace{0.1cm}\Bigg\vert \hspace{0.1cm} \tau >n\right] = \int_M h(y) \eta(\d y) \label{eq:conv}
\end{equation} \label{defqem}
holds for $\eta$-almost every $x\in M$.
\end{definicao}

\begin{mtheorem}\label{thm:QEM} Let $X_n$ be a  Markov chain on $E_M$ absorbed at $\partial$ satisfying Hypothesis \ref{(H)}. Let $\mu$ denote the unique quasi-stationary measure for $X_n$ and $\lambda$ its survival rate,  as in Theorem \ref{thm:QSM}. Let $f\in \mathcal C^0_+(M)$ be a non-negative continuous function such that $\mathcal P f = \lambda f$ and $m$ be the number of eigenvalues of $\mathcal P$ in the circle of radius $\lambda$, as defined in Theorem \ref{thm:Spec}. Assume additionally that either $m=1$ or $\rho(Z)=0$. Then $X_n$ admits a unique quasi-ergodic measure on $M$ given by $$\eta(\d x) = \frac{f(x) \mu(\d x)}{\int_M f(y) \mu(\d y)},$$
and equation \eqref{eq:conv} holds for every $x\in M\setminus Z$. Moreover,
\begin{enumerate}
    \item[(M1)] if $m=1,$ given for every $\nu \in \mathcal M_+(M),$ such that $\int f \d \nu >0$, there exist constants $K(\nu)$, $\alpha >0,$ such that
   \begin{align}
     \left\|\mathbb P_{\nu}\left[X_n\in\cdot \mid \tau >n \right]  - \mu \right\|_{TV} \leq K(\nu)e^{-\alpha n}, \ \fa \ n\in \mathbb N;  \label{expconvh}
   \end{align}
    \item[(M2)] if $\rho(Z) = 0,$  there exist open sets (on the induced topology of $M$) $C_{0}$, $C_{1}$, $\ldots$, $C_{m-1}$, such that $M\setminus Z = C_0 \sqcup C_1 \sqcup \ldots \sqcup C_{m-1},$ satisfying ${\{\mathcal P(\cdot,C_i)\neq 0\} = C_{i-1}},$   for all $i\in \{0,1,\ldots, m-1\},$ with $C_{-1} = C_{m-1}$. 
\end{enumerate}
\end{mtheorem}
Theorem \ref{thm:QEM} is proved in Section \ref{ProofC}.

 Observe that items (M1) and (M2) of Theorem \ref{thm:QEM} provide important information on the expected behaviour of $X_n$ before being absorbed into $\partial$.

In particular, if $m=1$, then for an $x\in M\setminus Z,$ we have that
$$\frac{\mathbb P\left[X_n \in A \mid X_0 =x\right]}{\mathbb P\left[X_n \in M \mid X_0 =x \right]} = \frac{\mathcal P^n(x,A)}{\mathcal P^n(x,M)} \to \mu(A),$$
exponentially fast, as $n\to\infty$. This limit means that the expected long-time behaviour of the noise realisations that stay in $M$ is described by the measure $\mu$. On the other hand, in case (M2), the process $X_n$ displays a cyclic behaviour.

\subsection{Applications}\label{sec:RDS}

In this subsection, we discuss some concrete applications of Theorems \ref{thm:QSM}, \ref{thm:Spec} and \ref{thm:QEM}. The main purpose of this subsection is to illustrate that the above theorems can be applied to a wide class of Markov chains. 

% The first example considers a Markov chain generated by a homeomorphism $T:\mathbb R\to\mathbb R$ and additive noise. 

\begin{exemplo}
\label{EXAMPLE:boundednoise}
  Consider $T: \mathbb R\to \mathbb R$ and $\varepsilon>0$, such that:
\begin{itemize}
\item $T$ is an increasing homeomorphism;
\item There exist unique $x_{-}$ and $x_{+}$ satisfying $T(x_-) + \varepsilon = x_-$ and ${T(x_+) - \varepsilon = x_+}$; and
\item $T(x) + \varepsilon \leq x$ for every $x<x_-$ and $x \leq T(x)-\varepsilon$ for every $x>x_+$.
\end{itemize}

Set $M = [x_-,x_+]$ and consider the Markov chain $X_{n+1}= T(X_n) + \omega_n$ where $\{\omega_i\}_{i\in\mathbb N_0}$ is an i.i.d. sequence of random variable uniformly distributed in $[-\varepsilon,\varepsilon]$. From this construction, it follows that $X_n$ is absorbed at $\partial = \mathbb R\setminus M$.

Given $A\in \mathscr{B}(M)$,
\begin{align*}
    \mathcal P(x,A) &= \frac{1}{2\varepsilon}\int_{-\varepsilon}^{\varepsilon} \mathbbm{1}_A (T(x)+\omega)\mathrm{d}\omega=\frac{1}{2|x_+ - x_-|\varepsilon}\int_{A} \mathbbm {1}_{M}( y - T(x)) \mathrm{d}y.
\end{align*}

It is readily verified that $X_n$ fulfils Hypothesis \ref{(H)} with $Z = \{x_-,x_+\}$, the only two points that immediately escape $M$ almost surely. Since $M\setminus Z$ is connected. Theorem \ref{thm:Spec} guarantees that $m=1$ and Theorems \ref{thm:QSM} and \ref{thm:QEM} provide the existence of a unique quasi-stationary and quasi-ergodic measure on $M$. Moreover, Theorem \ref{thm:QEM} also asserts that for every Borel measure $\nu$ on $M,$ such that $\supp (\nu) \not\subset Z,$
 $$ \left\|\mathbb P_{\nu}\left[X_n\in\cdot \mid \tau >n \right]  - \mu \right\|_{TV} \to 0\ \text{when }n\to \infty,$$
 exponentially fast.
\label{2x}

\end{exemplo}

Natural choices of $T$ in the above example are $T(x) = a x^{2n+1}$ for $n\in\mathbb N_0$, $a\in \mathbb R$ and $\varepsilon >0$ large enough. We note that in the case that $|T(x) - x|< 1$ for every $x\in M$, the existence of quasi-stationary and quasi-ergodic measures also follow from \cite[Proposition 7.4]{Champ2} and \cite[Corollary 6.5]{V2}. However, our results are applicable to a broader class of functions. In contrast, we would like to point out that the existence of quasi-stationary measures in the case where $M$ is not compact and $T(x) = a x$ was already investigated in \cite{Kolb1}.

In the next two examples, we explore natural extensions of Example \ref{EXAMPLE:boundednoise}. Unlike the previous scenario, which dealt with Markov chains generated by bounded additive perturbations of homeomorphisms in $\mathbb R$, Example \ref{x3} focuses on bounded geodesic-induced perturbations of the antipodal map in $\mathbb S^2$, and Example \ref{complex} studies absorbing Markov chains generated by perturbations of quadratic polynomials conditioned upon staying in a neighbourhood of its Julia set.

\begin{exemplo} \label{x3}
Let $E = \mathbb S^2$ be the unit sphere in $\mathbb R^3$ and let $M = \{(x,y,z)\in \mathbb S^2; x\geq 0\}$ be a hemisphere of $E$. Fix $\varepsilon >0$ and consider the Markov chain
$$X_{n+1} = \begin{cases}
    -Q(\varphi_n,\theta_n) X_n &, \text{if } X_n\in M\\
    X_n &, \text{if } X_n\in E\setminus M
    \end{cases},  $$ where  $$ Q(\varphi_n,\theta_n):=
 \begin{bmatrix}
\cos(\varphi_n) & -\sin(\varphi_n) & 0 \\
\sin(\varphi_n) & \cos(\varphi_n) & 0 \\
0 & 0 & 1
\end{bmatrix} \begin{bmatrix}
1 &0& 0 \\
0 & \cos(\theta_n) & -\sin(\theta_n) \\
0 & \sin(\theta_n) & \cos(\theta_n)
\end{bmatrix}, $$
and $\{(\varphi_n,\theta_n)\}_{i\in\mathbb N_0}$ is an i.i.d. sequence of random variable uniformly distributed in $[-\pi/2,\pi/2]\times [-\varepsilon,\varepsilon]$. By construction, this process effectively stops as soon as $X_n \in \partial = E\setminus M$.

It is clear from the above construction that $Z= \{(1,0,0)\}$. For every $A\in\mathscr B(M),$ we have
\begin{align*}
    \mathcal P(x,B) &= \frac{1}{2 \pi \varepsilon  }\int_{-\pi/2}^{\pi/2} \int_{-\varepsilon}^{\varepsilon} \mathbbm{1}_A (-Q(\varphi,\theta)x)\mathrm{d}\varphi \d \theta.
\end{align*}

By employing the coordinate transformation $v=-Q(\varphi,\theta) x$ and making use of the rank theorem \cite[Theorem 4.1]{Lee2012IntroductionManifolds}, it follows that $X_n$ satisfies Hypothesis \ref{(H)} when considering $\rho$ as the volume form in $\mathbb S^2$. Furthermore, from the definition of $X_n$ and the aforementioned observations, it can be shown that $X_n$ satisfies Hypothesis \ref{(H)}.

Once again, since $M\setminus Z$ is connected, Theorem \ref{thm:Spec} ensures that $m=1$. Moreover, Theorems \ref{thm:QSM} and \ref{thm:QEM} imply that on $M$, $X_n$ admits a unique quasi-stationary measure which is fully supported, and also a unique quasi-ergodic measure. Furthermore, \eqref{expconvh} holds.

\end{exemplo}

\begin{exemplo}
\label{complex}
Given $c\in\mathbb C$, define $p_c:\mathbb C \to \mathbb C$ as $p_c(z) = z^2+ c$. Denote by $K_c := \{z: \|f^n(x)\| \not\to\infty\ \text{as }n\to\infty  \}$ and let $J_c := \partial K_c$ be the so-called \emph{Julia set}. $J_c$ is non-empty compact and totally invariant, i.e. $J_c = p_c(J_c) = p_c^{-1}(J_c)$ (\cite[Lemma 3.1.]{milnor}). 

Let $c\in \mathbb C$ be such that $J_c$ is connected, i.e. $c$ lies in the Mandelbrot set. Let $\delta>0$ and define $J^\delta_c =\{z\in J_c; \|z-c\|\leq \delta\}$.  Given $\varepsilon>0$, define the Markov chain $X_{n+1} = p_c(X_n) + \mathbbm 1_M(X_{n})\omega_n$ on $E:= \mathbb C$, where $\{\omega_i\}_{i\in\mathbb N_0}$ is an i.i.d. sequence of random variables uniformly distributed in $B_\varepsilon (0) :=\{z\in \mathbb C;\|z\|\leq \varepsilon\}$, and $M$ is the closure of the maximal subset of $\mathbb C$ such that: 
\begin{itemize}
    \item $J_c \subset M\subset J_c^\delta;$ and
    \item $X_n$ is transitive in the interior of $M$.
\end{itemize}Observe that $M$ is well-defined since $p_c:J_c\to J_c$ has a dense orbit \cite[Corollary 3.11.]{milnor}, moreover $J^{\delta_1}_c\subset M$ for $\delta_1>0$ small enough.

Given $A\in \mathscr B(M)$ we obtain that
\begin{align*}
    \mathcal P(x,A) &= \frac{1}{\pi\varepsilon^2}\int_{B_\varepsilon(0)} \mathbbm{1}_A (p_c(x)+\omega)\mathrm{d}\omega=\frac{1}{\pi\varepsilon^2} \int_{A} \mathbbm {1}_{B_0(\varepsilon)} (z-p_c(x)) \mathrm{d}z.
\end{align*}
The above equation and the definition of $M$ imply that $X_n$ satisfies \ref{(H)}. Since $p_c: J_c\to J_c$ is topologically exact \cite[Lemma 3.7.]{milnor}, we obtain that $X_n$ does not display cyclic behaviour in $M$, so that Theorem \ref{thm:QEM} implies $m=1$. Finally, it follows from Theorems \ref{thm:QSM} and \ref{thm:QEM} that, $X_n$ admit unique quasi-stationary measure and quasi-ergodic measures on $M$. As before, \eqref{expconvh} holds.

\end{exemplo}

In the last example, we show that the constant $m$ given by Theorem \ref{thm:Spec} can be strictly less than the number of connected components of $M\setminus Z$.
\begin{exemplo}\label{example22}
Let $E=\mathbb R$ and $M$ be a compact regular subset of $E$, i.e. $M = \overline{\mathrm{Int}(M)}.$  Consider $X_n$ be the Markov chain $X_{n+1} = X_n + \mathbbm 1_{M}(X_n) \omega_n$,   where $\{\omega_i\}_{i\in\mathbb N_0}$ is an i.i.d. sequence of random variables normally distributed with mean $0$ and variance $1$. Let $\mathbb P$ denote the standard Gaussian distribution.

Then, for every $A \in \mathscr{B}(M),$
\begin{align*}
    \mathcal P(x,A)&= \int_{\mathbb R} \mathbbm 1_A\left(x+ \omega\right)\mathbb{P}(\d \omega)=\frac{1}{\sqrt{2 \pi}}\int_{A} e^{-\frac{(y-x)^2}{2}}\d y
\end{align*}
and Hypothesis \ref{(H)} is fulfilled with $Z = \emptyset$. Since $X_n$ does not display cyclic behaviour in $M$, Theorem \ref{thm:QEM} implies that $m=1$. Therefore, by Theorems \ref{thm:QSM} and \ref{thm:QEM} we conclude that $X_n$ admits a unique quasi-stationary stationary measure which is fully supported on $M,$ and a unique quasi-ergodic measure.

Note that, in this case, $m=1$ while although $M$ may have an infinite number of connected components. 

\end{exemplo}

\section{Some direct consequences of Hypothesis \ref{(H)} \label{sec:P} } 
 The purpose of this section is to present three fundamental results that are extensively used throughout this paper. The first two are direct consequences of Hypothesis \ref{(H)}, while the third one is functional-analytic.

We begin by stating some properties of the map $\mathcal P$, which follow from standard arguments that can be found in the literature \cite{Tweedie,Ale}.

\begin{proposicao} \label{NaVi}
If $X_n$ fulfils Hypothesis \ref{(H)}, then the map $\mathcal P: L^\infty (M)\to L^\infty  (M)$, $\mathcal Pf= \mathbb E_x[f\circ X_1 ]=\int_M f(y) \mathcal P(x,\d y)$ has the following properties:
\begin{enumerate}
    
    \item[(a)] For all $n\in\mathbb N$ and $f\in L^\infty(M),$ $(\mathcal P)^{n}(f) = \mathbb E_x[f \circ X_n];$ 
    
    \item[(b)] If $f\in L^\infty(M),$ then $\mathcal P f \in \mathcal C^0(M);$
    
    \item[(c)] $\left.\mathcal P\right|_{\mathcal C^0(M)}:\left(\mathcal C^0(M),\|\cdot\|_\infty\right)\to \left(\mathcal C^0(M),\|\cdot\|_\infty\right)$ is {a} positive compact operator; and

    \item[(d)] If $z\in Z$ then $\mathcal P(z,M\setminus Z) = 0$.
\end{enumerate}
\end{proposicao}
 \begin{proof}
{Item $(a)$ follows from the Markov property of the process $X_n.$  To prove $(b),$ let $x\in M$ and $\varepsilon>0.$ By Hypothesis (H2) there exists  $\delta>0$ such that
 $d(x,z)<\delta$ then $\|g(x,\cdot) - g(z,\cdot)\|_1 <\varepsilon.$ Therefore, given $z\in M$, such that $d(x,z)<\delta$, we obtain
 \begin{align}
     \left|\mathcal P f(z) -\mathcal P f(x)\right|  &\leq \|f\|_{\infty}  \|g(x,\cdot) - g(z,\cdot)\|_1  <\varepsilon \|f\|_{\infty} \label{ub}
 \end{align}
 implying that $\mathcal P f \in \mathcal C^0(M).$
 
Item $(c)$ follows from inequality \eqref{ub} and the Arzelá-Ascoli Theorem \cite[Proposition 5.3.]{Ale}.

Finally, we prove $(d)$. Suppose by contradiction that there exists $z\in Z$, such that $\mathcal P(z,M\setminus Z) >0.$ Therefore, we can find a compact set $K\subset M\setminus Z,$ satisfying $\mathcal P(z, K) >0.$ In this way, given $k\in\mathbb N$, we obtain that
$$\mathcal P^k(z,M) = \int_{M}\mathcal P^{k-1}(x,M) \mathcal P(z,\d y)\geq \mathcal P(z,K) \min_{x\in K} \mathcal P^{k-1}(x,M)>0,$$
since, for every $k\in \mathbb N,$ $\mathcal P^{k-1}(\cdot, M)$ is a continuous function and $\mathcal P^k(x,M)>0$ for every $x\in K\subset M\setminus Z.$ This contradicts $z\in Z$.}
\end{proof}

\begin{proposicao}\label{matheustrick}
Let the absorbing Markov chain $X_n$ satisfy Hypothesis \ref{(H)}. Then the following assertions hold:
\begin{itemize}
    \item[(a)]  If there exists $x_0\in M\setminus Z,$ such that $\mathcal P(x_0,M)<1,$ then there exist $n_0\in\mathbb N$ and $\alpha \in (0,1)$ such that
    $$ \mathcal P^{n}(x,M)<\alpha^{\left\lfloor \frac{n}{n_0}\right\rfloor},\ {\   \text{for all}\ } x\in M.$$
    In particular, $\lim_{n\to\infty}\mathbb P_x\left[\tau >n\right]=\lim_{n\to\infty} \mathcal P^n(x,M) = 0,$ for all $x\in M,$
    and $r(\mathcal P)< 1.$
    \item[(b)] Let $\{A_i\}_{i\in \mathbb N} \subset \mathscr{B}(M)$ such that $\lim_{n\to \infty }\mathbbm 1_{A_n} = 0,$ $\rho\text{-}$almost surely. Then, $\lim_{i\to\infty }\|\mathcal P(x,A_i)\|_{\infty}=0.$
\end{itemize}
\end{proposicao}
\begin{proof}

 Since $\mathcal P(\cdot,M)$ is continuous, there exists an open neighbourhood $B$ of $x_0,$ such that ${\sup_{y\in B} \mathcal P(y,M)<1.} $

To prove $(a)$, note that given $x \in M \setminus Z$, by Hypothesis (H2), there exists $n_x = n(x,B),$ such that $\mathcal P^{n_x}(x,B) >0.$ Therefore
\begin{align*}
    \mathcal P^{n_x +1} (x,M)&= \int_B   \mathcal P(y,M) \mathcal P^{n_x}(x,\mathrm{d}y)  + \int_{M\setminus B}   \mathcal P(y,M) \mathcal P^{n_x}(x,\mathrm{d}y) \\
    &\leq \mathcal P^{n_x}(x,B) \sup_{y\in B}\mathcal P(y,M) + \mathcal P^{n_x}(x,M\setminus B)<1,
\end{align*}
which implies that $\mathcal P^{n_x +1} (x,M)<1$. Hence, given any $x\in M$, there exists $m_x= n_x+1$, such that $\mathcal P^{m_x}(x,M) <1$. Since $\mathcal P^{m_x}(\cdot,M)$ is continuous, there exists  an open neighborhood  $U_x$ of $x$ such that $\mathcal P^{m_x}(y,M)<1,$ for all   $y\in U_x.$

Thus  $M = \bigcup_{x\in M} U_x,$ and since $M$ is compact, there exist $x_1,x_2,\ldots,x_n\in M$ such that $M = U_{x_1} \cup U_{x_2} \cup \cdots \cup U_{x_n}.$ It is readily verified that taking $n_0 = m_{x_1}\cdot m_{x_2} \cdot\ldots\cdot m_{x_n}$ we have that $\mathcal P^{n_0}(y,M) <1.$  Hence, defining $\alpha = \sup_{x\in M} \mathcal P^{n_0 }(x,M)<1,$ it follows that
\begin{align*}
    \mathcal P^{n n_0} (x,M) &= \int_M \mathcal P^{n_0}(y,M) \mathcal P^{n_0(n-1)} (x, \mathrm{d} y) \leq \alpha \mathcal P^{(n-1)n_0}(x,M) \leq \ldots\leq \alpha^n.
\end{align*}

It remains to prove $(b).$ From the Arzelá-Ascoli,  $\{\mathcal P(\cdot, A_i)\}_{i\in\mathbb N}$ is pre-compact in $(\mathcal C^0(M),\|\cdot\|_\infty).$ (H1) implies that  for a fixed $x\in M,$ $\lim_{n\to\infty} \mathcal P(x,A_i)=0.$ Therefore, if $f$ is an accumulation point of ${\{\mathcal P(\cdot, A_i)\}_{i\in\mathbb N}}$ in $(\mathcal C^0(M),\|\cdot\|_\infty)$, $f(x) =0$ for every $x\in M.$ We can therefore conclude that $\mathcal P(\cdot,A_i)$ converges to $0$ in $\mathcal C^0(M).$
\end{proof}

The following proposition states a functional-analytic result that is extensively used throughout this paper.

\begin{proposicao}\label{mula}
Let $E$ be a Banach space, $\lambda$ a positive real number, and $T:E\to E$ a bounded linear operator such that $r(T) <\lambda.$ Then there exist constants $K>0$ and $\delta \in (0,\lambda)$ such that $\|T^n\| \leq K (\lambda-\delta)^n,$ for all $n\in\mathbb N$. 
\end{proposicao}

We do not include a proof of this result as it follows from standard arguments in the literature and the formula of the spectral radius $r(T) = \lim_{n\to \infty}\|T^n\|^{\frac{1}{n}}$ for bounded linear operators $T: E\to E.$

\section{Banach lattices \label{4}}

In this section, we introduce Banach lattices, which are essential for the proof of the main theorems of this paper. We show that the operator $\mathcal P$ is well-behaved from a Banach lattice point of view, allowing us to deduce important properties of the spectrum of the operator $\mathcal P$. We start this section with some basic definitions from Banach lattice theory and state some key results from this field, which we then apply to the operator $\mathcal P$.

{Given $(L,\leq)$ a partial ordered set and a set $B\subset L,$ we define, if they exists
$$\sup B = \min\{\ell\in L;\ b\leq \ell,{\   \text{for all}}\ b\in B\}\ \text{and }\inf B = \max\{\ell\in L;\ \ell \leq b ,{\   \text{for all}}\ b\in B\}. $$
We say that $L$ is a \emph{lattice}, if for every $f_1,f_2\in L,$
$$f_1\lor f_2 := \sup \{f_1,f_2\} ,\ f_1\land f_2:= \inf \{f_1,f_2\}  $$
exist. Additionally,  in the case that $L$ is a vector space and the lattice $(L,\leq)$ satisfies
$$f_1\leq f_2 \Rightarrow f_1+f_3\leq f_2+f_3,\ {\   \text{for all}}\  f_3\in L,\ \text{and} $$
$$f_1\leq f_2 \Rightarrow \alpha f_1 \leq \alpha f_2,\ {\   \text{for all}\ } \alpha>0, $$
then $(L,\leq)$ is called \emph{vector lattice}. Finally, if  $(L,\|\cdot\|)$ is a Banach space and the vector lattice  $(L,\leq )$ satisfies
$$|f_1|\leq |f_2|\Rightarrow \|f_1\|\leq \|f_2\|,$$
where $|f_1| := f_1\lor (-f_1),$ then the triple $(L,\leq,\|\cdot\|)$ is called a 
\emph{Banach lattice}. When the context is clear, we denote the Banach lattice $(L,\leq,\|\cdot\|)$ by $L$.

This paper uses two concepts from the Banach lattice theory: an \emph{ideal} of a Banach lattice, and an  \emph{irreducible operator} acting on a Banach lattice. A vector subspace $I \subset L$, is called an \emph{ideal} if, for every $f_1,f_2\in L$ such that $f_2\in I$ and $|f_1|\leq |f_2|,$ we have $f_1\in I$. A positive linear operator $T:L\to L$ is called \emph{irreducible} if,  $\{0\}$ and $L$ are the only $T$-invariant closed ideals of $T$.}

The following two results give us tools to analyse the spectrum of a compact positive irreducible operator. In section \ref{BL}, we show that these results apply to the operator $\mathcal P$ when restricted to a specific subspace of $\mathcal C^0(M)$. This procedure allows us to analyse the spectrum of the operator $\mathcal P$.

\begin{proposicao}[{\cite[Proposition~2.1.9~(iii)]{Ideal1}}]  \label{cone}
Let $M$ be a compact Hausdorff space. Consider the Banach lattice $\mathcal C^0(M)$. Then, $I$ is an ideal of $\mathcal C^0(M),$ if and only if $I = \left\{f\in \mathcal C^0(M);\ \left. f\right|_A = 0\right\},$
for some closed set $A$.
\end{proposicao}

\begin{definicao}
Let $L$ be a Banach lattice. We denote by $L_+ := \{f\in L;\ 0\leq f\}$, and write $$L^*_+ :=\{\varphi \in L^*; \ {\varphi}\ \text{is a continuous positive linear operator} \},$$
{where $L^*$ is the topological dual space of $L.$}

A point $f\in L$ is called \emph{quasi-interior} if $f\in L_+$ and, for every $\varphi$ $ \in L^*_+\setminus \{0\}$ we have $\varphi(f) > 0$.
\end{definicao}
% {
% \begin{remark}
% The concept of quasi-interior point will be only be used during the proof 
% \end{remark}}

The main results from the theory of Banach lattices used in the upcoming sections are summarised in the following theorem.

% \begin{teorema}[ {\cite[Propositions 5.2 and 5.3]{B1}} ]Let $L$ be a Banach lattice and suppose that $T$ is positive and $T^n$ is compact for some $n\in\mathbb N$. If $T$ is an irreducible operator, then $r(T)>0$ and $r(T)$ is an eigenvalue of $T$ of multiplicity one. Moreover, the eigenspace is spanned by $u$, a quasi-interior point. \label{A1}
% \end{teorema}

\begin{teorema}[{\cite[Proposition 5.2 and Theorem 5.3]{B1}}] Let $L$ be a Banach lattice and let $T$ be an irreducible operator. If $T^k$ is compact for some $k\in \mathbb N$, then $r(T)>0$, and if $\lambda_1,\lambda_2,\ldots,\lambda_m$ are the different eigenvalues of $T$ satisfying $|\lambda_j|=r=r(T)$, for $j=1,\ldots,m,$ every $\lambda_j$ is a root of the equation $\lambda^m -r^m =0$. All these eigenvalues have algebraic multiplicity one and the spectrum of $T$ is invariant under a rotation of the complex plane under the angle $2\pi/m$, multiplicities included. \label{A2}  Moreover, the eigenspace $\mathrm{ker}(T-r(T))$ is spanned by a quasi-interior point $u\in L_+$.
\end{teorema}

\subsection{Banach lattice properties of \texorpdfstring{$\mathcal P$}{P}\label{BL}}

In this section we exploit the Banach lattice of $(\mathcal C^0(M), \leq ,\|\cdot\|_\infty)$ as the stochastic Koopman operator $\mathcal P$ acts from $\mathcal C^0(M)$ to itself, i.e. $\mathcal P:\mathcal C^0(M)\to\mathcal C^0(M).$ It is shown in Sections \ref{QSM1} and \ref{6} that the existence and uniqueness of quasi-stationary and quasi-ergodic measures are closely related to the spectrum of $\mathcal P$, $\mathrm{Spec}(\mathcal P)$. To our advantage, the tools described above allow us to give a detailed characterisation of $\mathrm{Spec}(\mathcal P)$, which we present in the following theorem.

% In this section, we exploit the Banach lattice property of the operator $\mathcal P$ to analyse its spectrum. Throughout this section we consider $\mathrm{Dom}(\mathcal P) =\mathcal C^0(M)$ and the Banach lattice $(\mathcal C^0(M),\leq , \|\cdot\|_\infty)$  induced by the natural Banach lattice structure of $\mathcal C^0(M)$.

% The following theorem describes the spectrum of $\mathcal P$. It is shown later that the existence and uniqueness of quasi-stationary and quasi-ergodic measures are closely related to the spectrum of $\mathcal P$.

\begin{teorema} \label{3o}

Let $X_n$ be a  Markov chain on $E_M$ absorbed at $\partial$ satisfying Hypothesis \ref{(H)} and consider the operator $\mathcal P :\mathcal C^0(M)\to \mathcal C^0(M).$ Then,  defining $\lambda := r(\mathcal P) >0,$ there exists $m\in \mathbb N$ such that $\left\{\lambda e^{\frac{2\pi i j}{m} }\right\}_{j=0}^{m-1}$ are the unique eigenvalues of modulus equal to $\lambda$, and each one has algebraic multiplicity one.

Moreover, $\mathrm{dim}\left(\mathrm{ker}\left(\mathcal P - \lambda e^{\frac{2\pi i j}{m}}\right)\right)=1,\ {\   \text{for all}}\ j \in \{0,1,\ldots,m-1\},$
and there exists $f\in \mathcal C_+^0(M),$ such that $f(x) =0$ if and only if $x\in Z.$

\end{teorema}

\begin{proof}

From Proposition \ref{NaVi} $(d)$ and Hypothesis \ref{(H)}, we obtain that $Z = \{x\in M, \mathcal P(x,M\setminus Z) = 0\}$. Since $x\mapsto \mathcal P(x,M\setminus Z)$ is continuous, then the set $ Z$ is compact.  Moreover, defining the set $\mathcal C_Z := \{f\in\mathcal C^0(M);\ f(z) =0,$ for all $ z\in Z\},$ it is clear that $\mathcal C_Z$ is a closed invariant subspace of $\mathcal C^0(M)$ and therefore $\mathcal C_Z$ admits a  Banach lattice structure induced by $\mathcal C^0(M).$ From Hypothesis \ref{(H)} and Propositions \ref{NaVi} and \ref{cone} it is readily verified that $\left.\mathcal P\right|_{\mathcal C_Z}:\mathcal C_Z\to \mathcal C_Z$ is an irreducible compact positive operator. Therefore, applying Theorem \ref{A2} to the operator $\left.\mathcal P\right|_{\mathcal{C}_Z}$, we have that $\left.\mathcal P\right|_{\mathcal C_Z}$ fulfils all the conditions stated in Theorem \ref{3o}. The proof is finished noticing that $\mathrm{Spec}(\mathcal P) = \mathrm{Spec}\left(\left.\mathcal P\right|_{\mathcal C_Z}\right),$ since $f$ if is an eigenvector of $\mathcal P$ then $f(z) = 0$ for every $z\in Z.$
\end{proof}

\section{{Existence and Uniqueness of a Quasi-stationary measures}\label{QSM1}}

In this section, we show that Hypothesis (H) is sufficient to obtain existence and uniqueness of a quasi-stationary measure for $X_n$ on $M$.

Recall that $\mathcal M(M):=\left\{\mu; \ \mu\ \text{is a Borel signed measure on }M\right\}$ has a Banach space structure when endowed with the total variation norm. It is well known, from the Riesz–Markov–Kakutani representation theorem \cite[Theorem 6.19]{rudin}, that $(\mathcal C^0(M),\|\cdot\|_{\infty})^*=\left(\mathcal M\left(M\right),\|\cdot\|_{TV}\right).$ Thus we may identify any $\mu \in \mathcal M(M),$ with an element of $(\mathcal C^0(M),\|\cdot\|_{\infty})^*,$ by  $\mu(f) := \int f(x) \mu(\d x),$ for every $f\in\mathcal C^0(M).$

In order to prove the existence and uniqueness of a quasi-stationary measure for $X_n$ we study the spectrum of the operator $\mathcal{L}: \mathcal{M}(M)\to \mathcal M (M)$, $ \mathcal L\mu := \int_M \mathcal P(x,\cdot) \mu(\d x).$ However, it is first necessary to linearly extend the action of $\mathcal P$ and $\mathcal L$ to their complex domains
$$\mathcal C^0(M,\mathbb C) := \mathcal C^0(M) \oplus i \mathcal C^0(M)\ \text{and }\mathcal M(M,\mathbb C):= \mathcal M(M)\oplus i \mathcal M(M). $$

\begin{definicao}[\cite{brezis}]
Let $E$ be a  Banach space. Then, the \emph{scalar product for the duality $E^*$, $E$} is the bilinear form  $\langle \cdot, \cdot \rangle: E^* \times E \to \mathbb C$, defined by $\langle \varphi , v\rangle := \varphi(v). $

Given a bounded linear operator $T: E\to  E$, we denote $T^* : E^* \to E^*$ as the linear operator $T^*\varphi (v)=  \varphi(T v)$, for all $v\in E$.
\end{definicao}

 Using this definition, it follows that $\mathcal P^* = \mathcal L$.  The following lemma shows a connection between the spectrum of the operators $\mathcal P$ and $\mathcal L.$ The proof follows from standard functional analytical techniques.

\begin{lema}
The operators $\mathcal P$ and $\mathcal L$ have the same eigenvalues and $$\mathrm{dim}(\mathrm{ker} (\mathcal L - {\beta})) = \mathrm{dim}(\mathrm{ker} (\mathcal P - {\beta} )),\ {\   \text{for all}} \ {\beta} \in \mathrm{Spec}(\mathcal P) =\mathrm{Spec}(\mathcal L). $$ 

Moreover, if $f_0\in \mathcal C^0_+(M)$ is an eigenfunction of $\mathcal P$ with eigenvalue ${\lambda}=r(\mathcal P)=r(\mathcal L),$ then there exists an eigenmeasure $\mu_{f_0}\in \mathcal M_+(M)$ of $\mathcal L$ with respect to the  eigenvalue ${\lambda}$ such that $\mu_{f_0}(f_0)=1. $
% where $$\mathcal M_+(M):=\left\{\mu \in \mathcal M(M); \ \mu (f)\geq 0,\ \text{for every }f\in\mathcal C_+^0(M)\right\}.$$
\label{autovalores}
\end{lema}

We can finally state and prove the main result of this section.
\begin{teorema}
If a Markov chain $X_n$ satisfies Hypotehsis \ref{(H)}, then $X_n$ admits a unique quasi-stationary {measure} $\mu$ on $M,$ and $\mathrm{supp} (\mu)$ satisfies \eqref{support}. \label{QSM}
\end{teorema}
\begin{proof}

We divide the proof into three steps.
\begin{step}{1}
{\it We show that if $\mu \in \mathcal M_+(M)$ is an eigenmeasure of $\mathcal L,$ then $$\mathrm{supp} \left(\mu\right)\supset \bigcup_{k\in \mathbb N}\bigcup_{x\in M\setminus Z} \supp(\mathcal P^k(x,\d y))\supset M\setminus Z.$$}\label{supp}
\end{step}
 
First, observe that there exists $x_0\in \supp(\mu)\cap M\setminus Z.$ In fact, if $\mu(Z)= 1$, define $Z_k := \{x\in Z; \mathcal P^k(x,M)=0\}.$ It is clear that $Z = \bigcup_{k\in\mathbb N} Z_k$. Therefore for every $k\in \mathbb N$
\begin{align*}
    1 = \mu(Z) &= \frac{1}{\lambda^k}\int_{Z} \mathcal P^k(z,Z)\mu(\d z) = \int_{Z\setminus Z^k} \frac{1}{\lambda^k}\mathcal P^k(z,Z)\mu(\d z)\\
    &\leq \sup_{z\in Z}\left|\frac{1}{\lambda^k} \mathcal P^k(z,Z)\right| \mu(Z\setminus Z_k) \xrightarrow[]{k\to \infty} 0.
\end{align*}
The last inequality holds since $r(\mathcal P) = \lambda$ and $\mu(Z_k) \to \mu(Z)$ as $k\to\infty.$

Suppose by contradiction that there exists an open set $A\subset M\setminus Z$, such that $\mu(A) =0$.   Therefore, for every open neighbourhood $U$ of $x$, $\mu(U)>0.$ 

Using Hypothesis (H2), there exists $n>0$ such that $\mathcal P^n(x_0,A)>0.$ Since $\mathcal P^n(\cdot,A)$  is continuous there exists an open neighbourhood $B$ of $x$ such that $\mathcal P^{n}(y,A) >\mathcal P^n (x_0, A)/2>0,$ for all $y\in B.$

On the other hand,
\begin{align*}
    0 &= \mu(A)= \frac{1}{\lambda^n} \int_M \mathcal P^n(y,A) \mu(\d y)\geq \frac{1}{\lambda^n} \frac{\mathcal P^n(x_0,A)}{2}\mu(B)>0,
\end{align*}
which is a contradiction since $x_0\in\supp(\mu)$ and $x_0\in B.$ {Thus} $\supp(\mu) = M\setminus Z.$ The proof is finished by observing that $$\mu(\d y) \geq \frac{1}{\lambda^k}\int_{M\setminus Z} \mathcal P^k(x,\d y) \mu(\d x)\ \text{for every }k\in\mathbb N.$$

\begin{step}{2}\label{step2}
{\it We show that the operator $\mathcal L$ admits a unique eigenmeasure that lies in $\mathcal M_+(M).$} 
\end{step}
  Define $\lambda := r(\mathcal P)=  r(\mathcal L)$. Combining Theorem \ref{3o}$ and Lemma  \ref{autovalores},$ we have that
\begin{align}
    \mathrm{dim}(\mathrm{ker}(\mathcal L - \lambda )) =\mathrm{dim}(\mathrm{ker}(\mathcal P - \lambda )) = 1\label{escalada}
\end{align}
and there exists a probability measure $\mu \in\mathcal M_+(M),$ such that $\mathcal L\mu = \lambda \mu.$ 

Suppose by contradiction that there exists a Borel probability measure $\nu$ on $M$ linearly independent from $\mu$ such that $\mathcal L \nu = \lambda_0 \nu$. From Theorem \ref{A2}, the linear operator $\mathcal P$ admits a eigenfunction $f\in \mathcal C^0_+(M)$, such that $f(x)>0$ for every $x\in M\setminus Z$  with eigenvalue $\lambda = r(\mathcal P)=r(\mathcal L)$. Therefore, $0 < \nu(f) \leq  \|f\|_{\infty} \nu(M) < \infty.$ On the other hand,
    \begin{align*}
    \int_M f(x) \nu(\mathrm{d}x)&=\frac{1}{\lambda_0}\langle \mathcal L\nu, f \rangle =\frac{1}{\lambda_0}\langle \nu,\mathcal Pf  \rangle =\frac{\lambda}{\lambda_0}\int_{M}  f(x) \nu (\mathrm{d} x). 
        \end{align*}
Therefore $\lambda_0 = \lambda,$ contradiction \eqref{escalada}.  This concludes Step 2.

\begin{step}{3}  {\it We show that the  Markov chain $X_n$ admits a unique quasi-stationary measure $\mu,$ and $\mathrm{supp}(\mu)$ satisfies \eqref{support}.}
\end{step}
Let $\mu$ be the unique probability eigenmeasure of $\mathcal L$, given by Step \ref{step2}.  We claim that $\mu$ is a quasi-stationary measure. Note that, for every $A\in{\mathscr B}(M)$ and $n\in\mathbb N,$
$$\mathbb P_{\mu}\left[X_n \in A\mid \tau >n\right] = \frac{\int_M\mathcal P^n(x,A)\mu(\d x)}{\int_M\mathcal P^n(x,M)\mu(\d x)} = \frac{\mathcal L^n\mu(A)}{\mathcal L^n\mu(M)} = \frac{\lambda^n \mu (A)}{\lambda^n \mu(M)} = \mu(A), $$
showing that $\mu$ is a quasi-stationary measure for $X_n$.

Reciprocally, if $\nu$ is a quasi-stationary measure for $X,$ then $\nu\in \mathcal M_+(M)$. Therefore, defining $\lambda_0 = \int_M \mathcal P(x,M)\nu(\d x),$ we obtain $\int_M\mathcal P(x,\cdot)\nu(\d x) = \lambda_0 \nu(\cdot).$

Hence, there is a $1$-$1$ correspondence between probability eigenmeasures of $\mathcal L$ lying in $\mathcal M_+(M)$ and quasi-stationary measures of $X_n$.

This concludes the proof of the theorem.
\end{proof}

We close this section with the proof of Theorem \ref{thm:QSM}.

\begin{proof}[Proof of Theorem \ref{thm:QSM}]
Note that from Theorem \ref{QSM}, $X_n$ admits a unique quasi-stationary measure $\mu$ and \eqref{support} holds. We first prove part  $(a).$ Since $\mathcal P(x,M)=1$ for all $x\in M$, the constant function $x\mapsto 1$ is an eigenfunction of $\mathcal P$. Therefore, $ r(\mathcal P) = r(\mathcal L) =1$ and $\mu$ corresponds to a stationary measure. In case $(b)$, since there exists $x_0\in M\setminus Z$ such that $\mathcal P(x_0,M)<1,$ Proposition \ref{matheustrick}  $(a)$ guarantees that $\lim_{n\to\infty} \mathcal P^n(y,M) = 0,\ {\   \text{for all}} \ y\in M, $
implying that $\mu$ is the unique quasi-stationary measure for $X_n$ with survival rate $\lambda <1.$
\end{proof}

% \section{Proof of Theorem \ref{thm:QSM}}

\section{Existence of a Quasi-Ergodic Measure \label{6}}

The proof of the existence of a quasi-ergodic measure is much more intricate than that of a quasi-stationary measure. Our technique is inspired by the seminal reference \cite{QEM}, which focuses on time-continuous absorbing Harris' chains, and \cite{Esperanca}, where time-discrete Markov chains were considered in the context of moving boundaries.  In this section, we extend the results presented in \cite{QEM} to include discrete-time Markov processes, which may not necessarily be aperiodic.

As in \cite{Esperanca}, our method depends on the number of eigenvalues of $\mathcal P$ in the circle ${r(\mathcal P)\mathbb S^1 \subset \mathbb C.}$ For this reason, we define the following quantity.

\begin{notacao}
Given a Markov chain $X_n$ satisfying Hypothesis \ref{(H)}, from Theorem \ref{3o} it follows that the number of eigenvalues in the peripheral spectrum of $\mathcal P$, i.e. $\mathrm{Spec}(\mathcal P) \cap r(\mathcal P) \mathbb S^1$, is finite. For the remaining of this section we denote $\lambda := r(\mathcal P) = r(\mathcal L)$ and $m :=\#\{\alpha\in\mathrm{Spec}(\mathcal P);\ \|\alpha\|=\lambda\}$.\label{m}
\end{notacao}

\begin{proposicao}\label{niko}
Let $X_n$ be a Markov chain on $E_M$ absorbed at  $\partial$ satisfying Hypothesis \ref{(H)}, and $\lambda = r(\mathcal P) = r(\mathcal L)$. Consider $f_0$, $f_1$, $\ldots$, $f_{m-1} \in \mathcal C^0(M,\mathbb C)$ and $\mu_0$, $\mu_1$, $\ldots$, $\mu_{m-1} \in \mathcal M(M,\mathbb C )$ such that $\mathcal P f_j = \lambda e^{\frac{2i\pi j}{m}} f_j$ and $\mathcal L \mu_j = \lambda e^{\frac{2i\pi j}{m}} \mu_j,$  for all $j\in\{0,1,\ldots,m-1\}.$  Then $\langle \mu_{j} ,f_{r}  \rangle = 0, \text{if  } j\neq r.$
In particular, normalise $f_i$ with respect to $\mu_i$ such that $\langle \mu_j , f_k\rangle = \delta_{j k },$ where $\delta$ is the Kronecker delta.
\end{proposicao}
\begin{proof}
Observe that
\begin{align*}
 \langle \mu_j, f_r \rangle &= \frac{1}{\lambda e^{\frac{2\pi i j}{m}}} \langle \mathcal L  \mu_j , f_r\rangle =\frac{1}{\lambda e^{\frac{2\pi i j}{m}}} \langle  \mu_j ,\mathcal P f_r\rangle=\frac{\lambda e^{\frac{2\pi i r}{m}} }{\lambda e^{\frac{2\pi i j}{m}}} \langle  \mu_j , f_r\rangle= e^{\frac{2\pi i (r-j)}{m}} \langle  \mu_j , f_r\rangle.
\end{align*}
Given $j,r \in \{0,1,\ldots, m-1\}$, if $r\neq j$ we have that $e^{\frac{2\pi i (r-j)}{m}}\neq 1,$ and therefore  ${\langle \mu_j, f_r \rangle = 0.}$ In the following we argue that $\langle f_i,\mu_i\rangle >0$ for every $i\in \{0,1,\dots,m-1\}.$

From \cite[Theorem 8.4-5]{kr}, we can decompose  $\mathcal C^0(M,\mathbb C) =  W_0\oplus W,$ where $W_0 := \mathrm{span}_{\mathbb C} (f_0,f_1,\ldots,f_{m-1})$ and $W$ is a $\mathcal P$-invariant subspace of $\mathcal C^0(M,\mathbb C)$ satisfying $r\left(\left.\mathcal P\right|_W\right) < \lambda.$ Observe that the above decomposition implies that $\langle \mu_j , f_j\rangle \neq 0,$ some $j\in\{0,1,\ldots,m-1\},$ otherwise we would have that $\mu_j= 0.$ Redefining  $f_j$ as $f_j/\langle \mu_j,f_j\rangle,$ the proof is concluded. 
\end{proof}

Until the end of this section, we set $\{f_j\}_{j=0}^{m-1}\subset  \mathcal C^0(M,\mathbb C)$ and $\{\mu_j\}_{j=0}^{m-1}\subset \mathcal M(M,\mathbb C)$ as in Proposition \ref{niko} assuming that $\langle \mu_j , f_i \rangle = \delta_{ij}$ for every $i,j\in\{0,\ldots,m-1\}.$

Furthermore, from Theorem \ref{3o} and \cite[Theorem 8.4-5]{kr}, we have
\begin{eqnarray}
\mathcal C^0 (M,\mathbb C) = \mathrm{span}_{\mathbb C} (f_0,f_1,\ldots,f_{m-1}) \oplus W,\label{dec2}
\end{eqnarray}
where $W$ is $\mathcal P$-invariant subspace of $\mathbb C$ and $ r\left(\left.\mathcal P\right|_{W} \right) < \lambda$, and
\begin{eqnarray}
\mathcal M (M,\mathbb C) = \mathrm{span}_{\mathbb C} (\mu_0,\mu_1,\ldots,\mu_{m-1}) \oplus V, \label{dec1}
\end{eqnarray}
where $V$ is $\mathcal L$-invariant subspace of $\mathcal M(M,\mathbb C)$ and $ r\left(\left.\mathcal L\right|_{V} \right) < \lambda$.

Note that the  decompositions \eqref{dec2} and \eqref{dec1} imply that $\|\mathcal P^n\| = \mathcal O( \lambda^n)$ and  $\|\mathcal L^n\| = \mathcal O( \lambda^n)$. Indeed, writing $ \mathbbm 1_M = \alpha_0 f_0 + \ldots + \alpha_{m-1} f_{m-1} + w,$
where $\alpha_0,\ldots,\alpha_{m-1}\in\mathbb C$ and $w\in W.$ Since $\mathcal P$ is a positive operator,
\begin{align}
    \|\mathcal P^n\| &=  \left\|\sum_{i=0}^{m-1} \alpha_i \mathcal P^n f_i + \mathcal P^n w \right\|_{\infty}\leq \lambda^n  \sum_{i=0}^{m-1}  |\alpha_i| \|f_i\|_{\infty}+ \smallO (\lambda^n)= \mathcal O(\lambda^n). \label{lambdan}
\end{align}

Finally, since $\|\cL^n\| = \|\cP^n\|,$ for all  $n\in\mathbb N$ we find that $\|\cL^n\| = \mathcal O(\lambda^n)$. In the following lemma, we discuss the behaviour of the Dirac measures $\delta_x$ under the decomposition \eqref{dec1}.
% {Note that using decomposition \eqref{dec1},$ for every $x\in M,$ there exist a measure $\nu_x \in \mathcal V$ and real numbers $\alpha_i\in\mathbb R,
\begin{proposicao}[Decomposition of Dirac measures]\label{goodbye}
Let $x\in M$, then there exists $\nu_x \in V$ such that
\begin{eqnarray}
 \delta_x = f_0(x) \mu_0 + f_1(x)\mu_1 + \ldots + f_{m-1}(x)\mu_{m-1} + \nu_x.\label{MemoBolinha}
\end{eqnarray}
Moreover, the family of function $\{\nu_x\}_{x\in M}$ satisfies
\begin{align}
\sup_{x\in M}\|\nu_x\|_{TV} \leq 1 + \sum_{i=0}^{m-1}\|f_i\|_{\infty} \|\mu_i\|_{TV} < \infty, \text{and } \sup_{x\in M} \|\mathcal L^n \nu_x\|_{TV} = \smallO(\lambda^n).\label{goodbye1}
\end{align}

\end{proposicao}
\begin{proof}
By decomposition \eqref{dec1}, there exist $\alpha_0,$ $\ldots,$ $\alpha_{m-1}$ $\in$ $\mathbb C$, and $\nu_x$ $\in$ $V$ such that $\delta_x = \sum_{i=0}^{m-1}\alpha_i \mu_i + \nu_x.$ Since $\frac{1}{\lambda^n} \mathcal L^n \delta_x =\sum_{k=0}^{m-1} \alpha_k \mu_k e^\frac{2\pi i k  n}{m} + \frac{1}{\lambda^n} \mathcal L^n \nu_x, $ and  $\langle f_i,\mu_j \rangle = \delta_{i j}$ for every $j\in \{0,\ldots,m-1\}$, we {obtain} that
\begin{eqnarray}
 \left \langle \frac{1}{\lambda^n} \mathcal L^n\delta_x , f_j\right\rangle = \alpha_j e^{\frac{2\pi i j n}{m}} + \left \langle \frac{1}{\lambda^n} \mathcal L^{n}\nu_x, f_j\right\rangle. \label{Memoestrela}
\end{eqnarray}

On the other hand,
\begin{eqnarray}
 \left \langle \frac{1}{\lambda^n} \mathcal L^n\delta_x,f_j\right\rangle =\frac{1}{\lambda^n}\left \langle  \delta_x , \mathcal P^n (f_j)\right\rangle = \frac{\lambda^n e^{\frac{2\pi i j n}{m}}}{\lambda^n}\left \langle  \delta_x , f_j\right \rangle =  e^{\frac{2\pi i j n}{m}} f_j(x).\label{Memoestrela2}
\end{eqnarray}

{From \eqref{Memoestrela}$ and \eqref{Memoestrela2}$ it follows that}
$$f_j(x) = \left \langle \frac{1}{\lambda^{n m}} \mathcal L^{n m} \delta_x , f_j\right\rangle = \alpha_j  + \left \langle \frac{1}{\lambda^{n m } } \mathcal L^{n m} \nu_x , f_j\right\rangle ,\ {\   \text{for all}} \ n \in\mathbb N.$$ 
From Proposition \ref{mula} we obtain that $f_j (x) = \alpha_j + \lim_{n\to\infty}\left \langle \frac{1}{\lambda^{n m } } \mathcal L^{n m} \nu_x , f_j\right\rangle = \alpha_j. $

The last part of the proposition follows from the fact that
\begin{align*}
    \sup_{x\in M}\|\nu_x\|_{TV} &=\sup_{x\in M}\left\|\delta_x -\left( f_0(x)\mu_0 + \ldots + f_{m-1}(x) \mu_j\right)\right\|_{TV}\\
    &\leq \sup_{x\in M}\|\delta_x\|_{TV} + \sum_{j=0}^{m-1}\|f_j\|_{\infty} \|\mu_j\|_{TV}=1 + \sum_{j=0}^{m-1}\|f_j\|_{\infty} \|\mu_j\|_{TV}<\infty,
\end{align*}
and $ \frac{1}{\lambda^n}\sup_{x\in M}\left\|\mathcal L^n \nu_x\right\|_{TV} \leq \frac{1}{\lambda^n} \left\|\left.\mathcal L^n\right|_{V}\right\|_{TV} \sup_{y\in M}\|\nu_{y}\|\longrightarrow 0,\ \text{when }n\to\infty,$
due to Proposition \ref{mula}.
\end{proof}

% \begin{remark} Using a similar argument, it is possible to prove that given a measure $\sigma\in\mathcal M_+(M),$ there exists $\nu_\sigma \in V,$ such that
% $$\sigma(\d y) = \int_{M} f_0(x)\sigma(\d x) \mu_0(\d y) + \ldots +\int_M f_{m-1}(x) \rho(\d x) \mu_{m-1}(\d y) + \nu_\sigma(\d y). $$
% \end{remark}

The following lemma lies at the basis of our proof for the existence of a quasi-ergodic measure. Since the proof of this result is long and technical, we include it in Appendix \ref{app}.

\begin{lema}
Let $x\in M\setminus Z$, and $h:M \to \mathbb R$ a bounded measurable function. Then, for every $n\in\mathbb N$.
$$\mathbb E_x\left[\sum_{k=0}^{n-1} h(X_k) \mathbbm 1_M (X_n)\right] = n \lambda^n \sum_{\ell = 0}^{m-1}  e^{\frac{2\pi i n \ell }{m}  }f_\ell  (x) \langle \mu_\ell  ,h\cdot f_\ell  \rangle \mu_\ell  (M)\ + \smallO (n\lambda^n). $$\label{Lobo}
\end{lema}

In the following two subsections, we analyse the cases when the number of eigenvalues of the operator $\mathcal P$ in $\lambda \mathbb S^1$, $m$, is either one or greater than one. The first case is much simpler compared to the second one. Indeed, if $m =1$, the operator $\mathcal P$ has the spectral gap property, simplifying the situation. Meanwhile, in the case $m>1,$ the process $X_n$ admits a cyclic property, which requires a more sophisticated analysis.

\subsection{Analysis of the case \texorpdfstring{$m=1$}{m=1}\label{ProofB}} In this section, we use Lemma \ref{Lobo} in order to prove item (M1) of Theorem \ref{thm:QEM}.

\begin{teorema}\label{m1}
Let $X_n$ be a Markov chain satisfying Hypothesis \ref{(H)}, and assume that $m =1.$ Let $\mu$ be the unique quasi-stationary measure on $M$, given by Theorem \ref{thm:QSM}, and $f$ be the unique function $\mathcal C^0_+(M)$ satisfying $\mathcal P f= \lambda f$ and $\int_M f(x) \mu(\d x) = 1$. Then, the measure $\eta(\d x) = f(x) \mu(\d x)$ is a quasi-ergodic measure for $X_n$ on the set $ M$ and \eqref{eq:conv} holds for every $x\in M\setminus Z$. 

Moreover, given $\nu \in \mathcal M_+(M),$ such that $\int_M f(y) \nu(\d y) >0,$ there exist $K(\nu) >0$ and $\alpha>0,$ such that $\left\|\mathbb P_\nu\left[X_n \in \cdot \mid \tau > n\right] - \mu \right\|_{TV}< K(\nu)e^{-\alpha n},$  for all $n\in\mathbb N.$

\end{teorema}
\begin{proof}
Let $x\in M$ and $h\in{\mathcal F_b}(M).$ Recall from Definition \ref{defqem} that we need to show that
$$\lim_{n\to\infty}\mathbb E_x\left[\frac{1}{n} \sum_{i=0}^{n-1} h\circ X_i \hspace{0.1cm}\Bigg\vert \hspace{0.1cm} \tau >n\right] = \int_M h(y) \eta(\d y), $$
for $\eta$-almost every $x\in M$. Since $\eta(Z)=0$, let us consider $x\in M\setminus Z.$

From Lemma \ref{Lobo} we obtain that
\begin{align}
    \mathbb E_x\left[\sum_{k=0}^{n-1} h(X_k) \mathbbm 1_M (X_n)\right] &= n\lambda^n f(x) \int_M h(y) f(y) \mu(\d y) + \smallO(n\lambda^n).\label{Memocarinha}
\end{align}

On the other hand, from Proposition \ref{goodbye}, there exists $\nu_x \in V$ such that $\delta_x = f(x) \mu + \nu_x$ and $r\left(\left.\mathcal L\right|_V\right)<\lambda.$ Since $ \mathcal L^n(\nu_x)(M)  =  \smallO{(\lambda^n)},$ by Proposition \ref{mula}, {we obtain}
\begin{align}
 \mathcal P^n(x,M) & =\left\langle \mathcal L^n \delta_x, \mathbbm 1_M\right\rangle = \lambda^n f(x) \mu(M) + \mathcal L^n \nu_x (M)  =\lambda^n f(x)  + \smallO{(\lambda^n)}.\label{memo3}
\end{align}

Hence, from \eqref{Memocarinha} and \eqref{memo3}
\begin{align*}
\mathbb E_x\left[\frac{1}{n} \sum_{i=0}^{n-1} h\circ X_i \hspace{0.1cm}\Bigg\vert \hspace{0.1cm} \tau >n\right]&=  \frac{1}{n}\frac{\displaystyle\sum_{i=0}^{n-1} \mathbb E_x \left[h(X_i)\mathbbm 1_m\right]}{\mathcal P^n(x,M)}=\frac{1}{n}\frac{n \lambda^n f(x) \mu( h f) + \smallO{(n\lambda^n)} }{\lambda^n f(x) + \smallO{(\lambda^n)}}\\
    &\xrightarrow{n\to\infty} \mu(h  f) \mu(\d y) = \int_M h(y)\eta(\d y).
\end{align*}

Now, we prove the second part of the theorem. Let $\nu\in\mathcal M_+(M)$ be such that $\nu(f) =\int_M f(y) \nu(\d y) >0.$ Given any measurable set $A$, by Proposition \ref{goodbye}
\begin{align}
    \mathcal L^n(\nu)(A) &= \int_{M} \mathcal L^n\delta_x(A)\nu(\d x)=\lambda^n \mu(A) \nu(f) + \int_M \mathcal L^n\nu_x(A) \nu(\d x) \label{memo4}. 
\end{align}

From Proposition \ref{mula}, there exist $\widetilde{K}>0$ and $\delta\in (0,\lambda),$ such that $ \left\|\mathcal L^{n}\right|_{V}\| <  \widetilde{K} \left( \lambda - \delta \right)^n,$ for every $n\in \mathbb N_0$
and therefore we can define $$\alpha := \left|\log \left(\frac{\lambda - \delta }{\lambda} \right) \right|,\ \text{and } K(\nu) := \widetilde{K} \frac{\nu(M)}{\nu(f)}\left(\sup_{x\in M}\|\nu_x\|_{TV}\right) (1+\|\mu\|_{TV})<\infty.$$
Thus, from \eqref{memo4}, we obtain that for every $n\in\mathbb N,$
\begin{align*}
    \left\|\mathbb P_\nu\left[X_n \in \cdot \mid \tau > n\right] - \mu \right\|_{TV} =& \left\|\frac{\int_M \mathcal P^n(x,\cdot) \nu(\d x)}{\int_M \mathcal P^n(x,M) \nu(\d x) } - \mu\right\|_{TV}\\
    =& \left\|\frac{\lambda^n  \nu(f) \mu + \int_M \mathcal L^n \nu_x (\cdot) \nu(\d x)}{\lambda^n \nu(f) + \int_M \mathcal L^n \nu_x (M) } - \mu\right\|_{TV}\\
    \leq &\frac{1}{\lambda^n \nu(f) }\left\| \int_M \mathcal L^n \nu_x (\cdot) \nu(\d x) - \mu \int_M \mathcal L^n \nu_x (M) \nu(\d x) \right\|_{TV}\\
    \leq & \frac{\nu(M)}{\nu(f)}\left(\frac{\|\left.\mathcal L^n\right|_V\|}{\lambda^n}\right)\left(\sup_{x\in M}\|\nu_x\|_{TV}\right) \left(1+\|\mu\|_{TV}\right)\\
    &\leq K(\nu) e^{-\alpha n}.
\end{align*}
\end{proof}
% {
% \begin{remark}
% If we choose a measure $\nu \in\mathcal M_+(M),$ such that $\int f(x)\nu(\d x)>0,$ it is also possible to prove, with a similar argument, that
% $$\lim_{n\to\infty}\mathbb E_\nu\left[\frac{1}{n} \sum_{i=0}^{n-1} h\circ X_i \hspace{0.1cm}\Bigg\vert \hspace{0.1cm} \tau >n\right] = \int_M h(y) \eta(\d y),\   {\   \text{for all}} \ h\in {\mathcal F_b}(M). $$

% \end{remark}}

At this stage, a significant portion of Theorems \ref{thm:Spec} and \ref{thm:QEM} has already been proven. Theorem \ref{thm:Spec} still requires us to show that $m$ is less or equal to the number of connected components of $M\setminus Z,$ and Theorem \ref{thm:QEM} requires the proof of condition (M2). We will address these remaining aspects in the upcoming subsection.

\subsection{Analysis of the case \texorpdfstring{$m>1$}{m>1}\label{ProofC}} In this section we show that $m>1$ implies that the Markov chain $X_n$ admits a cyclic behaviour and a quasi-ergodic measure on $M$. Also, in order to prove the second part of Theorem \ref{thm:QEM}, we always assume that $\rho(Z) = 0$.

{To conduct the desired analysis we study linear transformations $\mathcal P^m$ and $\mathcal L^m$.  From  Proposition \ref{niko}, it is clear that $r(\mathcal P^m) = r(\mathcal L^m) = \lambda^m$ and  $\mathrm{Spec}\left(\mathcal P^m\right)\cap \lambda^m \mathbb S^1 =\mathrm{Spec}\left(\mathcal L^m\right)\cap \lambda^m \mathbb S^1 = \{\lambda^m\}.$} Moreover,it also follows that $\mathrm{dim}\left(\mathrm{ker}(\mathcal P^m - \lambda^m)\right) $ $=$ $\mathrm{dim}\left(\mathrm{ker}(\mathcal L^m - \lambda^m )\right)$ $=$ $m.$

Throughout this section, we denote the \emph{support of a function} $f\colon M\to\mathbb R,$ as  $\mathrm{supp}(f) = \{x\in M;\ f(x)\neq 0\}$. 

In the following proposition, we study the eigenfunctions of $\mathcal P^m$ associated to the eigenvalue $\lambda^m$. A consequence of the proposition below is that 
$$ m \leq \#\{\text{number of connected components of }M\setminus Z\}. $$

\begin{proposicao}\label{nerdola}
There exist eigenfunctions $g_0,\ldots,g_{m-1} \in \mathcal C^0_+(M)$ of the operator $\mathcal P^m$ such that $\int_{M} g_j \d \mu =1$ for every $j\in\{0,1,\ldots,m-1\}$ and $\s_{\mathbb C} (\{g_i\}_{i=0}^{m-1}) = \mathrm{ker}( \mathcal P^m - \lambda^m ). $

Moreover, the eigenfunctions $g_0,$ $g_1,$ $\ldots$, $g_{m-1}$  can be chosen to have disjoint support, i.e., with $C_i := \mathrm{supp}(g_i)$ for all $i\in\{0,\ldots,m-1\}, $ satisfying $ C_i\cap C_j = \emptyset$ for all $i\neq j.$

Furthermore, the family of sets $\{C_i\}_{i=0}^{m-1}$ satisfies $ M\setminus Z = C_0 \sqcup C_1 \sqcup \ldots \sqcup C_{m-1}.$ In particular, $m\leq \#\{\mbox{connected components of }M\setminus Z\}.$
\end{proposicao}

\begin{proof}
Observe that since $\lambda^m\in \mathbb R$ and $\mathcal P(\mathcal C^0(M))\subset \mathcal C^0(M),$ it follows that if $f\in \mathcal C^0(M,\mathbb C)$ satisfies $\mathcal P^m f = \lambda^m f,$ then $\mathcal P^m \mathrm{Re}(f) = \lambda^m \mathrm{Re} (f)$ and $\mathcal P^m \mathrm{Im}(f) = \lambda^m \mathrm{Im}(f)$.

Let $\mu$ be the unique quasi-stationary measure of $X_n$ given by Theorem \ref{thm:QSM}. Note that the operator $\mathcal P^m$ satisfies
$\int_M \frac{1}{\lambda^m}\mathcal P^m f(x)\mu(\d x) = \int_M f(x)\d \mu$ for every $f\in\mathcal C^0(M).$
By the same techniques of \cite[Theorems 3.1.1 and 3.1.3]{CLM} it follows that if $f\in\mathcal C^0(M)$ is an eigenfunction of $\mathcal P^m$ with eigenvalue $\lambda^m$, then $f^+(x) = \max\{0,f(x)\}$ and $f^-(x) = \max\{0,-f(x)\}$ are also eigenfunctions of $\mathcal P^m$ with $\lambda^m.$

We will divide the proof into two steps.

\begin{step}{1} \it 
We prove that if $h_1,h_2 \in \mathcal C^0_+(M)$ are eigenfunctions of $\mathcal P^m$ with  eigenvalue $\lambda^m$ such that  $G:=\{h_1>0\}\setminus \{h_2 >0\} \neq \emptyset,$ then  $ \mathbbm{1}_{G}h_1$ is an eigenfunction of $\mathcal P^m$ with eigenvalue $\lambda^m.$
\end{step}

By the observations at the beginning of the proof, for all $t \in \mathbb R_{\geq 0}$, $h_t = (h_1 - t h_2)^{+} $ is an eigenfunction for $\mathcal P^m,$ with eigenvalue $\lambda^m.$ 
% i.e. $\mathcal P^m \left((h_1 - t h_2)^{+}\right) = \lambda^m (h_1 - t h_2)^{+}$, for all $t>0.$

Note that, $(h_1 - s h_2)^+ \leq (h_1 - t h_2)^+, {\   \text{for all}} \ s>t. $ We claim, that {$\{(h_1 - t h_2)^+\}_{t\in \mathbb R_+}$ stabilises on $t$}, i.e. there exists $t_0\geq 0$, such that $(h_1 - t_0 h_2)^+ =(h_1 - t h_2)^+,$ for every $t>t_0. $

Suppose by contradiction that the above statement is false. Then, we can find $t_1<t_2<\ldots<t_{m+1},$ such that $$\mathrm{supp}\left((h_1 - t_{m+1}h_2)^+\right)\subsetneq \mathrm{supp}\left((h_1 - t_{m}h_2)^+\right) \subsetneq \ldots \subsetneq \mathrm{supp}\left((h_1 - t_{1}h_2)^+\right).$$

The above equation implies that $\{(h_1 - t_{1}h_2)^+, \ldots, (h_1 - t_{m+1}h_2)^+\}$ is a linearly independent set in $\mathcal C^0(M)$. This generates a contradiction since $\mathcal P^m$ admits only $m$ eigenfunctions with eigenvalue $\lambda^m$.

Therefore, $\{(h_1 - t h_2)^+\}_{t\in\mathbb R_{\geq 0}}$ stabilises at some $t_0$. Finally, $\lim_{t \to \infty} (h_1 - t h_2)^+ = \mathbbm{1}_{G} h_1 = (h_1 - t_0 h_2)^+, $
and $(h_1 - t_0 h_2)^+$ is an eigenfunction of $\mathcal P^m.$

\begin{step}{2} \it We complete the proof of the theorem.

\end{step}
By Step 1 and the observations in the first two paragraphs, from linear combinations of the positive and negative parts of $\{f_i\}_{i=0}^{m-1}$ one can construct normalised  non-negative eigenfunctions $g_0,\ldots,g_{m-1}$ $\in$  $\mathcal C^0_+(M)$ of $\mathcal P^m,$ such that $\s_{\mathbb C} (\{g_i\}_{i=0}^{m-1}) = \mathrm {ker}( \mathcal P^m - \lambda^m ),$ and the family of sets $\{C_i :=$ $\supp (g_i)\}_{i=0}^{m-1}$
fulfils $C_i\cap C_j = \emptyset\ {\   \text{for all}} \ i\neq j.$

To prove the last part of the proposition, let $f \in \mathcal C^0_+(M)$ be an eigenfunction of $\mathcal P$ with eigenvalue $\lambda.$ From Corollary \ref{3o},  $\supp(f) =M\setminus Z.$ Since $f\in \mathrm{ker}\left(\mathcal P^m - \lambda^m\right)$ $=$ $\s_{\mathbb C}\left(\{g_i\}_{i=0}^{m-1}\right),$
there exist $\alpha_0,\ldots,\alpha_{m-1}\geq 0,$ such that $f = \alpha_0 g_0 +\ldots + \alpha_{m-1} g_{m-1},$ implying that $M\setminus Z=C_0 \sqcup \ldots \sqcup C_{m-1}.$ {Since each $C_i$ is open and closed in the topology in induced by $M\setminus Z$, we have that $m\leq \#\{\text{connected components on  }M\setminus Z\}.$} 
\end{proof}

\begin{proof}[Proof of Theorem \ref{thm:Spec}] \ The proof follows directly from Theorem \ref{3o} and Proposition \ref{nerdola}.
\end{proof}

From now on, we consider$\{g_i\}_{i=0}^{m-1}\subset \mathcal C_+^0(M)$  as in Proposition \ref{nerdola} and $\{C_i:= \supp{(g_i)}\}_{i=0}^{m-1} $. We may decompose 
\begin{eqnarray}
\mathcal C^0(M) = \s_{\mathbb R}(g_0,g_1,\ldots,g_{m-1}) \oplus V_0,\label{dec}
\end{eqnarray}
where $r\left(\left.\mathcal P^m\right|_{V_0}\right)<\lambda^m.$ For convenience we denote $C_i = C_{i\ (\mathrm{mod}\ m)},$ for every $i\in\mathbb N_0.$

We proceed to address the cyclic property of the eigenvectors of $\mathcal P^m$ and $\mathcal L^m,$ when $m>1.$ The following four results show that we may choose suitable eigenfunctions and eigenmeasures of the operators $\mathcal P^m$ and $\mathcal L^m,$ {respectively}, so that these permute cyclically, by application of $\mathcal P$ and $\mathcal L$.

\begin{proposicao}

 Let $j\in \{0,1,\ldots,m-1\},$ then $\mathrm{supp}\left(\mathcal P^m(\cdot,A) \right) \subset C_j,$  for all  $A\in{\mathscr B}(C_j).$ Assume that $\rho(Z) = 0$ and let $\mu$ be the unique probability measure on $\mathcal M(M)$ such that $\mathcal L\mu =\lambda \mu$. Then, for every $i\in\{1,\ldots,m-1\},$ $\nu_i(\d x):= \mu(C_i\cap \d x)$ is an eigenmeasure of $\mathcal L^m$ with eigenvalue $\lambda^m.$

\label{carladias}
\end{proposicao}
\begin{proof}
% Let $n$ be a natural number. It is clear that $ g_j \leq \|g_j\|_{\infty} \mathbbm{1}_{C_j}$. Applying $\mathcal P^{nm}$ to the last equation, it follows that $0< \lambda^{mn}g_j(x) \leq \|g_j\|_{\infty} \mathcal P^{nm}(x,C_j),$ for every $x\in M$ and $n\in\mathbb N.$  Implying \eqref{star}.

Let $A\in\mathscr{B}(C_j)$. For every $n\in\mathbb N,$ let us define $A_n := A\cap \{g_j >1/n\}.$ Since $\mathcal P^m$ is a positive operator, it follow that for every $n\in\mathbb N$, $\mathcal P^m(x,A_n) \leq n \lambda^m  g_j(x)$. This implies that $\mathrm{supp}\left(\mathcal P^m(\cdot,A_m) \right) \subset C_j$ for every $n\in \mathbb N.$ Since $C_i$ is an open set and $\supp(\mathcal P(\cdot,A))\subset M\setminus Z,$ Proposition \ref{matheustrick} (b) completes the proof. 

For the last part of the proposition, observe that the above paragraph implies that for every $B \in \mathscr{B}(M)$ and $i\in\{1,\ldots,m-1\},$
$$\lambda^m \nu_i (B)=\lambda^m \mu (B\cap C_i) = \int_M \mathcal P^m (x,B\cap C_i)\mu(\d x) = \int_M \mathcal P^m (x,B)\nu_i(\d x). $$
In the last equality, we used that $\{\mathcal P^m(x,B\cap C_i)>0\}\cap C_j\} = \emptyset$ if $i\neq j$, and that $\mathcal P^m(x,B\cap Z) =0$ since $\rho(Z)=0.$ 
\end{proof}

From now on, we consider $\{\nu_i\}_{i=0}^{m-1}$ as in Theorem \ref{carladias}. The following two lemmas establish a combinatorial structure in the set of measures $\{\nu_i\}_{i=0}^{m-1}$ (and the set of functions $\{g_i\}_{i=0}^{m-1}$). We omit the proof of the first lemma since it follows directly from standard linear algebra techniques.
\begin{lema}
Let $A= [a_{ij}]_{i,j =1}^{m}$ be $m\times m$ matrix, such that the following properties hold:
\begin{enumerate}
    \item[(i)] $A^m =\mathrm{Id};$
    \item[(ii)] $\displaystyle\sum_{j=1}^m a_{ij} \leq 1, \ {\   \text{for all}}\ i\in \{1,\ldots,m\}; and$
    \item[(iii)] $a_{ij}\geq 0, {\   \text{for all}}\  i,j\in\{1,\ldots,m\}.$
\end{enumerate}

 Then, the matrix $A$ permutes the canonical basis. In other words, for every $k\in\{1,\ldots,m\}$ there exists $s\in\{1,\ldots,m\}$, such that $ Ae_k = e_s.$\label{permutacao}
\end{lema}

\begin{lema}
Let $\{g_i\}_{i=0}^{m-1}\subset \mathcal C^0_+(M)$ as in \eqref{dec} and $\{\nu_i\}_{i=0}^{m-1}\subset \mathcal M_+(M)$ as in Theorem \ref{carladias}. Then,
there exists a cyclic permutation  $\sigma:\{0,1,\ldots,m-1\}\to\{0,1,\ldots,m-1\}$ of order $m$ such that for every $i\in\{0,1,\ldots,m-1\},$ $\frac{1}{\lambda}\mathcal L \nu_i = \nu_{\sigma(i)}$ and $\frac{1}{\lambda}\mathcal P g_i = g_{\sigma^{-1}(i)}.$ Moreover $\supp \,\mathcal P(\cdot ,C_{i}) = C_{\sigma^{-1}(i)}$ for every $i\in\{0,1,\ldots,m-1\}.$

\label{ciclying}
\end{lema}
\begin{proof}

We divide the proof into four steps.

\begin{step}{1}
\it There exists a permutation $\sigma$ on $\{0,1,\ldots,m-1\}$ such for every $i\in\{0,1,\ldots,m-1\},$ $\frac{1}{\lambda}\mathcal P g_i = g_{\sigma(i)}$ and $\supp \, \mathcal P(\cdot, C_i)= C_{\sigma^{-1}(i)}.$\label{estepe1}
\end{step}

Given $i\in\{0,1,\ldots,m-1\},$ note that $\mathcal P^{m}\left(\mathcal P g_i \right) =  \mathcal P(\mathcal P^{m} g_i ) = \lambda^m \mathcal P g_i,$ implying that $\mathcal Pg_i$ is an eigenmeasure of $\mathcal P^m$ with eigenvalue $\lambda^m$. Since $\mathcal P$ is a positive operator, $\mathcal P g_i\in \mathcal C_+^0(M).$ This implies that, for every {$i\in\{0,1,\ldots,m-1\}$} there exist $\alpha_{i0}, \ldots,\alpha_{i{m-1}}\in \mathbb R_+,$ such that $\frac{1}{\lambda}\mathcal P g_i = \sum_{j=0}^{m-1}\alpha_{ij} g_j. $

% Note that, since
% \begin{align*}
%     \mathrm{span}_{\mathbb C}(\nu_0,\nu_1,\ldots,\nu_{m-1})&= \bigoplus_{j=0}^{m-1} \mathrm{ker}\left(\cP-\lambda e^{\frac{2\pi j i }{m}}\right),  
% \end{align*}
% we have $\frac{1}{\lambda}\|\mathcal L \nu_i\|_{TV} \leq \left\|\frac{1}{\lambda} \left.\cL\right|_{\mathrm{span}_{\mathbb C}(\nu_0,\nu_1,\ldots,\nu_{m-1})} \right\|    \|\nu_i\|_{TV} \leq 1$ for all $i\in\{0,1,\ldots,m-1\}.$
Observe that
\begin{align*}
    1 &= \int g_i(x) \ \mu(\d x) = \int \frac{1}{\lambda} \mathcal P g_i(x) \ \mu(\d x) \\ & = \sum_{j=0}^{m-1} \alpha_{ij}  \int_{M} g_j(x) \mu(\d x) = \sum_{j=0}^{m-1} \alpha_{ij}
\end{align*} for all $i\in\{0,1,\dots,m-1\}$.  Defining the $m\times m$ matrix $A := \left[\alpha_{ij}\right]_{i,j=0}^{m-1},$ we obtain that $A^m = \mathrm{Id}.$

By Lemma \ref{permutacao}, the matrix $A$ permutes the canonical basis, and therefore, for every $i\in\{0,1,\ldots,m-1\},$ there exists $j_i\in\{0,1,\ldots,m-1\},$ such that 
$\frac{1}{\lambda}\mathcal P g_i  = g_{j_i}.$

Defining {$\sigma:\{0,1,\ldots,m-1\}\to \{0,1,\ldots,m-1\},$} such that $\sigma^{-1}(i) = j_i,$ we conclude the first step and $\mathrm{supp}\left( \mathcal P g_{j}\right) \subset C_{\sigma^{-1}(j)}.$

\begin{step}{2}
\it {$\sigma$ is a cyclic permutation of order $m$}.\label{estepe2}
\end{step}
{
Suppose that the permutation $\sigma$ admits a $k$-subcycle  $\widetilde{\sigma}.$  Without loss of generality assume that $\widetilde{\sigma} = (0,\widetilde{\sigma}(0),\ldots, \widetilde{\sigma}^{k-1}(0) ).$ Defining $\widetilde{f} = \frac{1}{k}\sum_{i=0}^{k-1}g_{\widetilde\sigma^{-i}(0)}$}, we have that $\mathcal P\widetilde{f} = \lambda \widetilde{f},$
implying that $\widetilde{f} \in \ker (\mathcal P -\lambda \mathrm{Id}).$ Since ${\mathrm{dim} \ker (\mathcal P -\lambda \mathrm{Id}) =1}$ we obtain that $f$ and $\widetilde f$ are linearly dependent, therefore $k=m$.

\begin{step}{3}\it We conclude the proof of the lemma. \label{es}
\end{step}

Let $\sigma$ be the $m$-cycle constructed in Steps \ref{estepe1} and \ref{estepe2}, such that
\begin{eqnarray}
  \frac{1}{\lambda} \mathcal P g_i = g_{\sigma^{-1}(i)},\ \text{for every }i\in\{0,1,\ldots,m-1\}.\label{permL}
\end{eqnarray}

It remains to show that $ \frac{1}{\lambda} \mathcal L \nu_i = \nu_{\sigma(i)},\ \text{for every }i\in\{0,1,\ldots,m-1\}.$ As in Step 1, for every $j\in\mathbb \{0,1,\ldots,m-1\}$, $\mathcal L \nu_j $ is an eigenfunction of $\mathcal L^m.$ Therefore, there exists $\beta_{j 0},\ldots,\beta_{j m-1}>0,$ such that
$\frac{1}{\lambda}\mathcal L \nu_j= \sum_{i=0}^{m-1}\beta_{ji} \nu_{i}. $
 From equation \eqref{permL} and duality, we obtain
$$  \frac{1}{\lambda}\langle \mathcal L \nu_{i}, g_j\rangle = \frac{1}{\lambda}\langle \nu_{i},\mathcal P g_j \rangle =\langle  \nu_{i}, g_{\sigma^{-1}(j)} \rangle=  \delta_{\sigma(i)  j} \langle \nu_{i},g_{\sigma^{-1}(j)} \rangle,$$ for every $i,j\in\{0,1,\ldots,m-1\}.$ Since $\frac{1}{\lambda}\mathcal L \nu_j= \sum_{i=0}^{m-1}\beta_{ji} \nu_{i}$ and $\langle \nu_i, g_j\rangle = \delta_{i j}$, we obtain that  $ \mathcal L v_{i} = \lambda  \nu_{\sigma(i)}$.

% \begin{step}{4}\it Proof of part $(b).$
% \end{step}

% Note that for every $k\in\{0,\ldots,m-1\},$
% \begin{align*}
% 1= \langle \nu_0, g_0\rangle&= \left\langle \frac{1}{\lambda^{m-k}} \cL^{m-k} \nu_0 , \frac{1}{\lambda^k}\mathcal P^k g_0 \right\rangle 
% = \left\langle \nu_{\sigma^{m-k}(0)}, g_{\sigma^{m-k}(0)}\right\rangle.
% \end{align*}

\end{proof}

% \begin{remark}\label{nb2}
The previous lemma establishes that $\mathrm{supp}(\mathcal P(\cdot,C_i)) = C_{\sigma^{-1}(i)},$ for all $i\in\{0,1,\ldots, m-1\},$ and $\mathrm{supp}(\mathcal P^k(\cdot, C_i))$ $\cap$ $\mathrm{supp}(\mathcal P^s(\cdot,C_i)) = \emptyset$ for all $k,s\in \mathbb N$ with $k\neq s$ $(\mathrm{mod}\ m).$ From Lemma \ref{ciclying}, we label the components $C_i$ such that

\begin{align}
    \mathcal P(x , C_i) & \begin{cases}
=0,\ \text{if} \ x \not\in C_{i-1}\\
>0,\ \text{if} \ x \in C_{i-1}
\end{cases},\label{j1}
\end{align} 
and for every $i\in \mathbb N_0,$
\begin{align}
     \frac{1}{\lambda}\mathcal Pg_i &= g_{i-1 \ (\mathrm{mod}\ m)} \ \text{and }\frac{1}{\lambda}\mathcal L \nu_i = \frac{1}{\lambda}\int_M \mathcal P(x,\cdot)\nu_i(\d x) = \nu_{i+1\ (\mathrm{mod}\ m)}. \label{j2}
\end{align}

We proceed to characterise the eigenvectors of the operators $\mathcal P$ and $\mathcal L$ and show that their eigenvalues lie in $\lambda \mathbb S^1\subset \mathbb C$.
\begin{teorema} Let $X_n$ be a  Markov chain on $E_M$ absorbed at $\partial$ satisfying (H) and $\rho(Z) = 0$. Let $\{g_i\}_{i=0}^{m-1}\subset \mathcal C^0_+(M)$  and $\{\nu_i\}_{i=0}^{m-1}\subset \mathcal M_+(M)$ be as in Theorem \ref{ciclying} and in \eqref{j1}$-\eqref{j2}$. Then, for every $j \in \{0,1,\ldots, m -1\}$, the measure
$$\mu_j = \frac{1}{m}\sum_{k=0}^{m-1} e^{\frac{-2\pi i k j}{m}} \nu_k, \text{ and the function }f_j = \sum_{k=0}^{m-1} e^{\frac{2 \pi i k j}{m}} g_k, $$
satisfy
\begin{align}
    \mathcal L\mu_j = \lambda e^{\frac{2\pi i j}{m}}\mu_j\ \mbox{
and }\mathcal Pf_j =  \lambda e^{\frac{2\pi i j}{m}}f_j. \label{feedback}
\end{align}

Moreover, the measures $\eta_{j} (\d x) :=  f_j (x) \mu_j(\mathrm{d}x)$ do not depend on $j\in\{0,1,\ldots,m-1\}$.\label{Dores1}
\end{teorema}
\begin{proof}

First, we verify that \eqref{feedback} holds. Consider $j\in\{0,1,\ldots,m-1\}$. From Theorem \ref{ciclying} and re-ordering of $\{\nu_i\}_{i=0}^{m-1}$ and $\{g_i\}_{i=0}^{m-1}$ we have that
\begin{align*}
    \mathcal L\mu_j& = \mathcal L\left(\frac{1}{m}\sum_{k=0}^{m-1} e^{\frac{-2\pi i k j}{m}} \nu_k\right) = \frac{1}{m}\sum_{k=0}^{m-1} \lambda e^{\frac{-2\pi i k j}{m}} \nu_{k+1} = \lambda e^{\frac{2\pi i j}{m}}\mu_j 
\end{align*}
and
\begin{align*}
    \mathcal P f_j& = \mathcal P\left(\sum_{k=0}^{m-1} e^{\frac{2\pi i k j}{m}} g_k\right) =\sum_{k=0}^{m-1} \lambda e^{\frac{2\pi i k j}{m}} g_{k-1} = \lambda e^{\frac{2\pi i j}{m}}f_j.
\end{align*}

From Lemma  \ref{ciclying} $(b)$, we have for every $j\in\{0,1,\ldots,m-1\}$ and $A\in\mathscr B(M)$
\begin{align*}
    \eta_j(A) &= \int_{A} f_j(x) \mu_j(\mathrm{d} x)=\frac{1}{m}\int_{A} \left(\sum_{k=0}^{m-1}e^{\frac{2\pi i k j} {m} }  g_k(x)\right)\left( \sum_{s=0}^{m-1} e^{\frac{-2\pi i j s}{m}}\nu_s(\mathrm{d} x)\right)\\
    &=\frac{1}{m}\sum_{k=0}^{m-1} \int_{A}  g_k(x) \nu_k(\mathrm{d} x) + \frac{1}{m}\sum_{k\neq s}e^{\frac{2\pi i j (k-s)}{m}}\underbrace{\int_{A}  g_k(x) \nu_s(\mathrm{d} x)}_\text{=0} \\
    &=\int_{A}  f_0(x) \mu_0 (\mathrm{d}x) = \eta_0(A),
\end{align*}
which completes the proof.
\end{proof}
% \begin{remark}
% In particular, Theorem \ref{Dores1} establishes that
% $\displaystyle \mu_0 = \frac{1}{m}\sum_{j=0}^{m-1}\nu_j,$
% is the unique quasi-stationary measure for $X_n,$ given by Theorem \ref{thm:QSM}.
% \end{remark}

\begin{corolario}
In the context of Theorem \ref{Dores1}, if the Markov chain $X_n$ satisfies (H) and $\rho(Z) =0$, then for every $x\in M\setminus Z$  and $h\in{\mathcal F_b}(M),$ 
\begin{align}
\mathbb E_x\left[\sum_{k=0}^{n-1} h(X_k) \mathbbm 1_M (X_n)\right] = n \lambda^n \langle \mu_0 , f_0 \cdot h \rangle \sum_{\ell = 0}^{m-1}  e^{\frac{2\pi i n \ell }{m}  }f_\ell  (x) \langle \mu_\ell   ,\mathbbm 1_M \rangle\ + \smallO (n\lambda^n). \label{asterisco}
\end{align}\label{rex2}
\end{corolario}
\begin{proof}
Combining Proposition \ref{Lobo},  Theorem \ref{Dores1}, and since $\langle \mu_\ell  ,h f_\ell  \rangle = \langle \mu_0,h f_0\rangle$ for every $\ell\in\{0,1,\ldots,$ $m-1\}$ we obtain \eqref{asterisco}.
\end{proof}

The following result establishes the existence of a quasi-ergodic measure for $X_n$ on $M$, in the case that $X_n$ satisfies Hypothesis \ref{(H)} and $\rho(Z)=0$.

\begin{teorema}\label{mch}
Let $X_n$ be a  Markov chain on $E_M$ absorbed at $\partial$ that satisfies Hypothesis (H) and $\rho(Z) = 0$. Then $X_n$ admits a quasi-ergodic measure $\eta$ on $M$ and \eqref{eq:conv} holds for every $x\in M\setminus Z.$

Moreover, $ \eta (\mathrm{d}x) =f_0(x)\mu_0(\d x), $ where $f_0 \in\mathcal C^0_+(M)$ and $\mu_0 \in \mathcal M_+(M)$, are such that $\mathcal P f_0 = \lambda f_0$, $\mathcal L \mu_0  = \lambda \mu_0,$ and $\langle \mu_0,f_0\rangle = 1$.

\end{teorema}
\begin{proof}
Let $h$ be a bounded measurable function. Recall from Definition \ref{defqem} that we need to show that
    $$\lim_{n\to\infty}\mathbb E_x\left[\frac{1}{n} \sum_{i=0}^{n-1} h\circ X_i \hspace{0.1cm}\Bigg\vert \hspace{0.1cm} \tau >n\right] = \int_M h(y) \eta(\d y), $$ 
for $\eta$-almost every $x\in M.$ Observe that $\eta(Z)=0$ and let $x\in X\setminus Z.$

From Proposition \ref{goodbye} we have
\begin{align}
    0 <\mathcal P^n (x,M)&=   \langle \mathcal L^n  \delta_x , \mathbbm 1_M \rangle=   \sum_{j=0}^{m-1} \lambda^n f_j (x) e^{\frac{2\pi i n j }{m}}\langle  \mu_j  , \mathbbm 1_M \rangle + \smallO(\lambda^n).\label{estrelinha}
\end{align}
Moreover, for every $n\in\mathbb N$
$$\mathbb E_x\left[ \frac{1}{n}\sum_{k=0}^{n-1} h(X_i)   \Bigg\vert \tau > n\right] = \frac{\displaystyle \mathbb E_x\left[\frac{1}{n}  \sum_{k=0}^{n-1} h(X_k)\mathbbm 1_M(X_n)\right]}{\mathcal P^n(x,M)}. $$
Applying Corollary \ref{rex2} for the numerator of the above equation, and \eqref{estrelinha} for the denominator, we obtain
\begin{align*}
    \mathbb E_x\left[ \frac{1}{n}\sum_{k=0}^{n-1} h(X_i)   \Bigg\vert \tau > n\right]& = \frac{\displaystyle \lambda^n \langle \mu_0,h\cdot f_0\rangle \sum_{\ell = 0}^{m-1}  e^{\frac{2\pi i n \ell }{m} } f_\ell  (x)  \langle \mu_\ell   ,\mathbbm 1_M \rangle + \smallO (n\lambda^n)}{\displaystyle \lambda^n \sum_{j=0}^{m-1} f_j (x) e^{\frac{2\pi i n j }{m}}\langle  \mu_j  , \mathbbm 1_M \rangle + \smallO(\lambda^n)}\\
    &=  \langle \mu_0,h\cdot f_0\rangle \frac{\displaystyle \sum_{\ell = 0}^{m-1}  e^{\frac{2\pi i n \ell }{m} } f_\ell  (x)  \langle \mu_\ell   ,\mathbbm 1_M \rangle + \smallO (1)_{n\to\infty}}{\displaystyle \sum_{j=0}^{m-1} e^{\frac{2\pi i n j }{m}}f_j(x) \langle  \mu_j  , \mathbbm 1_M \rangle + \smallO(1)_{n\to\infty}}.
\end{align*}
The proof concludes by taking limit  $n\to\infty$.
\end{proof}

We proceed to prove Theorem \ref{thm:QEM}.

\begin{proof}[Proof of Theorem \ref{thm:QEM}]
Observe that Theorems \ref{m1} and \ref{mch} imply that $\eta(\d x) = f(x) \mu(\d x)$ is a quasi-ergodic measure for $X_n$ on $M$ and \eqref{eq:conv} holds for every $x\in M\setminus Z.$  Moreover,  (M1) follows directly from Theorem \ref{m1} and (M2) is a consequence of Lemma \ref{ciclying}.

\end{proof}

\section*{Acknowledgments}

The authors gratefully thank Bernat Bassols-Cornudella for his valuable comments and suggestions. MC’s research has been supported by an Imperial College President’s PhD scholarship and FAPESP (process 2019/06873-2). MC and JL are also supported by the EPSRC Centre for Doctoral Training in Mathematics of Random Systems: Analysis, Modelling and Simulation (EP/S023925/1). JL gratefully acknowledges research support from IRCN, University of Tokyo and CAMB, Gulf University of Science and Technology, as well as from the London Mathematical Laboratory. GOM’s research has been supported by a PhD scholarship from CONACYT and the Department of Mathematics of Imperial College, and by the MATH+ postdoctoral program of the Berlin Mathematics Research Center.

\appendix
\section{Proof of Lemma \ref{Lobo}}
\label{app}
\label{appendixA}
This appendix contains a proof of Lemma \ref{Lobo}. The proof below is inspired by \cite[Proposition 4]{Esperanca}. Although the results in \cite{Esperanca} are focused on finite state spaces, it is possible to extend these results to our setting by making several adaptations. Since these are not straightforward, we present them here in detail.

The results of \cite{Esperanca} extend classical results by Darroch and Senata \cite{d}, where similar bounds are found for irreducible finite Markov chains.

\begin{proof}[{\it Proof of Lemma \ref{Lobo}}]

Fix $x\in  M\setminus Z$ and $h \in \mathcal F_b(M).$ We divide the proof into three steps.

\begin{step}{1}
{\it We prove that for every  $n\in\mathbb N,$} 
\begin{align}
    \sum_{k=0}^{n-1}\mathbb E_x \left[h(X_k) \mathbbm{1}_{M}(X_n)\right]=& \sum_{k=0}^{n-1} \sum_{\ell=0}^{m-1}\sum_{j=0}^{m-1} \lambda^{n} e^{\frac{2\pi i \ell k}{m} +\frac{2 \pi i (n-k)  j}{m}}f_\ell(x)\langle \mu_\ell, h f_j \rangle \langle \mu_j ,\mathbbm 1_M\rangle \label{Memosorriso}\\
     &+ \sum_{k=0}^{n-1}\left\langle \sum_{\ell=0}^{m-1}  \lambda^k e^{\frac{2\pi i \ell k}{m}} f_\ell  (x) \mu_\ell   , h \langle  \mathcal L^{n-k} \nu_{\cdot}, \mathbbm 1_M\rangle\right\rangle \nonumber\\
     &+\sum_{k=0}^{n-1} \left\langle  \mathcal L^k  \nu_x  , h \cdot \mathcal P^{n-k}\mathbbm 1_M\right\rangle.\nonumber
\end{align}
\end{step}

First, observe that
\begin{eqnarray}
 \mathbb E_x \left[ h(X_n) \mathbbm {1}_m\right] = \mathcal P^n h(x) = \langle \delta_x , \mathcal P^n h \rangle = \langle \mathcal L^n\delta_x,h\rangle. \label{Memobolinhabranca}
\end{eqnarray}

From the Markov property of $X_n$, for every $k,n\in\mathbb N,$ with $k\leq n,$ we have
\begin{align*}
    \mathbb E_{x} \left[h(X_k) \mathbbm 1_{M}(X_n)\right] &=\mathbb E_x \left[\mathbb E_{x} \left[h(X_k) \mathbbm 1_{M}(X_n)\mid \mathcal F_k\right]\right]= \mathbb E_x\left[h(X_k) \mathcal P^{n-k} (X_k , M)\right]\\&= \mathcal P^k ( h \cdot \mathcal P^{n - k}\mathbbm 1_M )(x).
\end{align*}
Using \eqref{Memobolinhabranca} and Proposition \ref{goodbye}  this implies that for every $y\in M$,
\begin{align*}
     h( y) \mathcal P^{n - k}(y , M)= h(y) \langle  \delta_y , \mathcal P^{n-k}\mathbbm 1_M\rangle =&  \sum_{j=0}^{m-1} \lambda^{n-k} e^{\frac{2 \pi i (n-k)  j}{m}} h(y)f_j(y)  \langle \mu_j,\mathbbm 1_M\rangle\\&+  h(y)\langle \mathcal L^{n-k}\nu_{y} , \mathbbm 1_M\rangle.
\end{align*}
{Recall} from Proposition \ref{goodbye} that $\sup_{y\in M}\|\mathcal L^n \nu_y\|_{TV} = \smallO(\lambda^n).$ {Hence,} 
\begin{align*}
     \mathbb E_x \left[ h(X_k) \mathbbm{1}_m\right]=& \mathcal P^k (h\cdot \mathcal P^{n-k}\mathbbm 1_M )(x)  = \left\langle \mathcal L^k\delta_x, h \cdot\mathcal P^{n-k}\mathbbm 1_M \right\rangle\\
     =&\left\langle \sum_{\ell=0}^{m-1}  \lambda^k e^{\frac{2\pi i \ell k}{m}} f_\ell  (x) \mu_\ell   , h\cdot \mathcal P^{n-k}\mathbbm 1_M \right\rangle + \left\langle  \mathcal L^k \nu_x , h \cdot \mathcal P^{n-k}\mathbbm 1_M\right\rangle \\
     =&\left\langle \sum_{\ell=0}^{m-1}  \lambda^k e^{\frac{2\pi i \ell k}{m}} f_\ell  (x) \mu_\ell   ,\sum_{j=0}^{m-1} \lambda^{n-k} e^{\frac{2 \pi i (n-k)  j}{m}} h f_j \langle \mu_j,\mathbbm 1_M\rangle +  h \langle \mathcal L^{n-k}\nu_y , \mathbbm 1_M\rangle  \right\rangle \\
     &+ \langle  \mathcal L^k \nu_x , h \cdot  \mathcal P^{n-k}\mathbbm 1_M\rangle \\
     =&\sum_{\ell=0}^{m-1}\sum_{j=0}^{m-1} \lambda^{n} e^{\frac{2\pi i \ell k}{m} +\frac{2 \pi i (n-k)  j}{m}}f_\ell  (x)\langle \mu_\ell  , h  f_j \rangle \langle \mu_j ,\mathbbm 1_M\rangle \\
     &+ \left\langle \sum_{\ell=0}^{m-1}  \lambda^k e^{\frac{2\pi i \ell k}{m}} f_\ell  (x) \mu_\ell   , h  \langle  \mathcal L^{n-k} \nu_{\cdot}, \mathbbm 1_M\rangle\right\rangle+ \langle  \mathcal L^k \nu_x , h\cdot \mathcal P^{n-k}\mathbbm 1_M\rangle. 
\end{align*}
Then, for each $n\in\mathbb N$, \eqref{Memosorriso} holds.  This completes  Step 1.

\begin{step}{2} {\it We prove that the following identity holds}
\begin{eqnarray}
 \label{Memotriste}
 \sum_{k=0}^{n-1}\left\langle \sum_{\ell=0}^{m-1}  \lambda^k e^{\frac{2\pi i \ell k}{m}} f_\ell  (x) \mu_\ell   , h \langle  \mathcal L^{n-k} \nu_{\cdot}, \mathbbm 1_M\rangle\right\rangle 
     +\sum_{k=0}^{n-1} \langle  \mathcal L^k \nu_x , h \mathcal P^{n-k}\mathbbm 1_M\rangle = \smallO (n \lambda^n).
\end{eqnarray}
\end{step}
 
% and equation $(\ref{lambdan})$ and  $\|\mathcal P^n\| = \mathcal O(\lambda^ n)$

Recall from Proposition \ref{goodbye} that  $\sup_{y\in M} \|\mathcal L^n \nu_y\| = \smallO(\lambda^ n)$. Hence, for all $\varepsilon>0$ there exists $n_0\in\mathbb N$ such that $n>n_0$ implies
\begin{align}
    \frac{1}{\lambda^n}\|\mathcal L^ n \nu_x\|_{TV} <\varepsilon.\label{eps}
\end{align}
On the other hand, recall from equation \eqref{lambdan} that $\|\mathcal P^n\| = \mathcal O(\lambda^n).$ Thus, there exists $K\geq 0$ such that $ \|\mathcal P^n\| \leq K \lambda^n$, for every $n\geq 0.$

We first discuss the second term in \eqref{Memotriste}. Note that for every $n>n_0+1$,
\begin{align*}
    \frac{1}{\lambda^ n n}\sum_{k=0}^{n-1} \langle  \mathcal L^k \nu_x , h\cdot \mathcal P^{n-k}\mathbbm 1_M\rangle \leq&  \frac{1}{\lambda^n n}\sum_{k=0}^{n-1} \|\mathcal L^k \nu_x\|_{TV} \|h\|_{\infty} \|\mathcal P^{n-k}\|\\
    % \leq&  \frac{K\|h\|_{\infty} }{\lambda^n n}\sum_{k=0}^{n-1} \|\mathcal L^k \nu_x\|_{TV} \lambda^{n-k}\\
    \leq&\frac{K\|h\|_{\infty} }{n}\sum_{k=0}^{n-1} \frac{\|\mathcal L^k \nu_x\|_{TV}}{\lambda^ k}\leq \frac{K\|h\|_{\infty} }{n}\left(\sum_{k=0}^{n_0-1} \frac{\|\mathcal L^k \nu_x\|_{TV}}{\lambda^ k} +\sum_{j=n_0}^{n-1} \varepsilon\right)\\
    &\longrightarrow K\|h\|_{\infty}\varepsilon,\ \text{as }n\to\infty.
\end{align*}
Hence, since $\varepsilon$ is arbitrary, 
 \begin{align}
     \sum_{k=0}^{n-1} \langle  \mathcal L^k \nu_x , h\cdot \mathcal P^{n-k}\mathbbm 1_M\rangle = \smallO (\lambda^n n). \label{stars} 
 \end{align}
 
 On the other hand, defining for each $n\in\mathbb N$
 $$I_n := \sum_{k=0}^{n-1}\left\langle \sum_{\ell=0}^{m-1}  \lambda^k e^{\frac{2\pi i \ell k}{m}} f_\ell  (x) \mu_\ell   , h\langle  \mathcal L^{n-k} \nu_{\cdot}, \mathbbm 1_M\rangle\right\rangle,$$
 we have
 \begin{align*}
     \frac{1}{n\lambda^ n} I_n \leq & \frac{1}{n\lambda^ n}\sum_{k=0}^{n-1} \sum_{\ell=0}^{m-1}  \lambda^k  \|f_\ell  \|_{\infty} \| \mu_\ell  \|_{TV} \|h\|_{\infty} \sup_{y\in M}\|\mathcal L^{n-k} \nu_{y}\|_{TV}\\
    %  \leq & C \frac{1}{n\lambda^n}\sum_{k=0}^{n-1}  \lambda^k \sup_{y\in M}\|\mathcal L^{n-k} \nu_{y}\|_{TV}\\
     =& C \frac{1}{n}\sum_{k=0}^{n-1}   \frac{\sup_{y\in M}\|\mathcal L^{n-k} \nu_{y}\|_{TV}}{\lambda^{n-k}}
      \end{align*}
 where $C:=m\max_{l \in\{0,\ldots,m-1\}}(\|h\|_\infty\|f_\ell  \|_{\infty} \|\mu_\ell  \|_{TV}).$
Hence, by \eqref{eps}
\begin{align*}  
\frac{1}{n\lambda^n}I_{n}     \leq& C\left( \frac{1}{n}\sum_{k=0}^{n_0-1}   \frac{\sup_{y\in M}\|\mathcal L^{k} \nu_{y}\|_{TV}}{\lambda^{k}} +\frac{1}{n}\sum_{k=n_0}^{n-1}   \frac{\sup_{y\in M}\|\mathcal L^{k} \nu_{y}\|_{TV}}{\lambda^{k}}\right)\\
     \leq & C\left( \frac{1}{n}\sum_{k=0}^{n_0-1}   \frac{\sup_{y\in M}\|\mathcal L^{k} \nu_{y}\|_{TV}}{\lambda^{k}} +\frac{1}{n}\sum_{k=n_0}^{n-1} \varepsilon\right)\\
     &\longrightarrow C \varepsilon,\ \text{when }n\to\infty,
 \end{align*}

 Once again, since $\varepsilon$ is arbitrary
 \begin{align}
\sum_{k=0}^{n-1}\left\langle \sum_{\ell=0}^{m-1}  \lambda^k e^{\frac{2\pi i \ell k}{m}} f_\ell  (x) \mu_\ell   , h \langle  \mathcal L^{n-k} \nu_{\cdot}, \mathbbm 1_M \rangle\right\rangle = \smallO(n\lambda^n).     \label{doublestars}
 \end{align} 
 
The sought identity is implies by \eqref{stars} and \eqref{doublestars}.
%  Proving that 
%  $$\sum_{k=0}^{n-1}\left\langle \sum_{\ell=0}^{m-1}  \lambda^k e^{\frac{2\pi i \ell k}{m}} f_\ell  (x) \mu_\ell   , h\cdot \langle  \mathcal L^{n-k} \nu_{\cdot}, \mathbbm 1_M\rangle\right\rangle 
%      +\sum_{k=0}^{n-1} \langle  \mathcal L^k \nu_x , h\cdot \mathcal P^{n-k}\mathbbm 1_M\rangle = \smallO (n \lambda^n). $$

\begin{step}{3}
{\it We prove that for every $n\in\mathbb N$.}
$$\mathbb E_x\left[\sum_{k=0}^{n-1} h(X_k) \mathbbm 1_M (X_n)\right] = n \lambda^n \sum_{\ell = 0}^{m-1}  e^{\frac{2\pi i n \ell }{m}  }f_\ell  (x) \langle \mu_\ell  ,h \cdot f_\ell  \rangle \mu_\ell  (M)\ + \smallO (n\lambda^n).$$
\end{step}
 
Note that by the previous two steps, $\sum_{k=0}^{n-1}\mathbb E_x \left[h(X_k) \mathbbm{1}_{M}(X_n)\right]= A
     + \smallO(n\lambda^ n),$
where
$$ A := \sum_{k=0}^{n-1} \sum_{\ell=0}^{m-1}\sum_{j=0}^{m-1} \lambda^{n} e^{\frac{2\pi i \ell k}{m} +\frac{2 \pi i (n-k)  j}{m}}f_\ell  (x)\langle \mu_\ell  , h f_j \rangle \langle \mu_j ,\mathbbm 1_M\rangle. $$

By exchanging the order of the sums, we obtain that
\begin{align*}
    %  A =& \sum_{\ell=0}^{m-1}\sum_{j=0}^{m-1} \left( \sum_{k=0}^{n-1} e^{\frac{2\pi i \ell k}{m} +\frac{2 \pi i (n-k)  j}{m}} \right)\lambda^{n} f_\ell  (x)\langle \mu_\ell  , h f_j \rangle \langle \mu_j ,\mathbbm 1_M\rangle\\
     A=&\sum_{\ell=0}^{m-1}\sum_{j=0}^{m-1} \left( \sum_{k=0}^{n-1} e^{\frac{2\pi i  n j}{m} +\frac{2 \pi i (\ell-j) k}{m}} \right)\lambda^{n} f_\ell  (x)\langle \mu_\ell  , h f_j \rangle \langle \mu_j ,\mathbbm 1_M\rangle.
\end{align*}
Splitting the above double sum into $\ell =j$ and $\ell \neq j$ we obtain
\begin{align*}
    A =& \sum_{\ell = 0}^{m-1} n   \lambda^n e^{\frac{2\pi i \ell n}{m}}f_\ell  (x) \langle \mu_\ell  ,h f_j\rangle \langle \mu_\ell   ,\mathbbm 1_M\rangle+\sum_{\ell \neq j}  \lambda^{n}e^{\frac{2 \pi i j }{m}} \left(\frac{e^{\frac{2\pi i \ell n}{m} } - e^{\frac{2\pi i j n}{m} } }{e^{\frac{2\pi i  \ell }{m} } - e^{\frac{2\pi i  j}{m} }}\right)f_\ell  (x) \langle \mu_\ell   , h f_j \rangle \langle \mu_j. \mathbbm 1_M\rangle
\end{align*} 
Since $e^{\frac{2 \pi i j }{m}} \left(\frac{e^{\frac{2\pi i \ell n}{m} } - e^{\frac{2\pi i j n}{m} } }{e^{\frac{2\pi i  \ell }{m} } - e^{\frac{2\pi i  j}{m} }}\right)$ is uniformly bounded in $n$ for $\ell,j\in\{0,1,\ldots,m-1\}$ and $\ell \neq j,$ we obtain that
$$\sum_{l \neq k}  \lambda^{n}e^{\frac{2 \pi i j }{m}} \left(\frac{e^{\frac{2\pi i \ell n}{m} } - e^{\frac{2\pi i j n}{m} } }{e^{\frac{2\pi i  \ell }{m} } - e^{\frac{2\pi i  j}{m} }}\right)f_\ell  (x) \langle \mu_\ell   , h f_j \rangle \langle \mu_j, \mathbbm 1_M \rangle  = \smallO (n\lambda^n). $$

The equation above implies
$$ \mathbb E_x\left[\sum_{k=1}^{n} h(X_k) \mathbbm 1_M (X_n)\right] = n \lambda^n \sum_{\ell = 0}^{m-1}  e^{\frac{2\pi i n \ell }{m}  }f_\ell  (x) \langle \mu_\ell  ,h f_\ell  \rangle \mu_\ell  (M)\ + \smallO (n\lambda^n).$$

This proves the lemma.
\end{proof}

%% if your bibliography is in bibtex format, uncomment commands:

\bibliographystyle{plain}
 
\bibliography{mybibfile.bib}

\begin{thebibliography}{10}

\bibitem{Araujo}
V.~Ara\'{u}jo.
\newblock Attractors and time averages for random maps.
\newblock {\em Ann. Inst. H. Poincar\'{e} C Anal. Non Lin\'{e}aire}, 17(3):307--369, 2000.

\bibitem{benaim2021}
M.~Bena\"{\i}m, N.~Champagnat, W.~O\c{c}afrain, and D.~Villemonais.
\newblock Transcritical bifurcation for the conditional distribution of a diffusion process.
\newblock {\em J. Theoret. Probab.}, 36(3):1555--1571, 2023.

\bibitem{QEM}
L.~A. Breyer and G.~O. Roberts.
\newblock A quasi-ergodic theorem for evanescent processes.
\newblock {\em Stochastic Process. Appl.}, 84(2):177--186, 1999.

\bibitem{brezis}
H.~Brezis.
\newblock {\em Functional analysis, Sobolev spaces and partial differential equations}.
\newblock Springer Science \& Business Media, 2010.

\bibitem{Kifer}
M.~Brin and Yu. Kifer.
\newblock Dynamics of {M}arkov chains and stable manifolds for random diffeomorphisms.
\newblock {\em Ergodic Theory Dynam. Systems}, 7(3):351--374, 1987.

\bibitem{collet2}
P.~Cattiaux, P.~Collet, A.~Lambert, S.~Mart\'{\i}nez, S.~M\'{e}l\'{e}ard, and J.~San~Mart\'{\i}n.
\newblock Quasi-stationary distributions and diffusion models in population dynamics.
\newblock {\em Ann. Probab.}, 37(5):1926--1969, 2009.

\bibitem{Champ}
N.~Champagnat and D.~Villemonais.
\newblock Exponential convergence to quasi-stationary distribution and {$Q$}-process.
\newblock {\em Probab. Theory Related Fields}, 164(1-2):243--283, 2016.

\bibitem{Champ2}
N.~Champagnat and D.~Villemonais.
\newblock General criteria for the study of quasi-stationarity.
\newblock {\em Electron. J. Probab.}, 28:Paper No. 22, 84, 2023.

\bibitem{QSB}
P.~Collet, S.~Mart{\'\i}nez, and J.~San~Mart{\'\i}n.
\newblock {\em Quasi-stationary distributions: {M}arkov chains, diffusions and dynamical systems}.
\newblock Springer Science \& Business Media, 2012.

\bibitem{martin1}
F.~Colonius and M.~Rasmussen.
\newblock Quasi-ergodic limits for finite absorbing {M}arkov chains.
\newblock {\em Linear Algebra Appl.}, 609:253--288, 2021.

\bibitem{d}
J.~N. Darroch and E.~Seneta.
\newblock On quasi-stationary distributions in absorbing discrete-time finite {M}arkov chains.
\newblock {\em J. Appl. Probability}, 2:88--100, 1965.

\bibitem{B2}
B.~de~Pagter.
\newblock Irreducible compact operators.
\newblock {\em Math. Z.}, 192(1):149--153, 1986.

\bibitem{Dennis1}
P.~Del~Moral and D.~Villemonais.
\newblock Exponential mixing properties for time inhomogeneous diffusion processes with killing.
\newblock {\em Bernoulli}, 24(2):1010--1032, 2018.

\bibitem{CarlangeloTransf}
M.~F. Demers, N.~Kiamari, and C.~Liverani.
\newblock {\em Transfer operators in hyperbolic dynamics---an introduction}.
\newblock 33 $^{\rm o}$ Col\'{o}quio Brasileiro de Matem\'{a}tica. Instituto Nacional de Matem\'{a}tica Pura e Aplicada (IMPA), Rio de Janeiro, 2021.

\bibitem{Op}
T.~Eisner, B.~Farkas, M.~Haase, and R.~Nagel.
\newblock {\em Operator theoretic aspects of ergodic theory}, volume 272 of {\em Graduate Texts in Mathematics}.
\newblock Springer Cham, 2015.

\bibitem{LExp}
M.~Engel, J.~S.~W. Lamb, and M.~Rasmussen.
\newblock Conditioned {L}yapunov exponents for random dynamical systems.
\newblock {\em Trans. Amer. Math. Soc.}, 372(9):6343--6370, 2019.

\bibitem{foguel}
S.~R. Foguel.
\newblock {\em The ergodic theory of {M}arkov processes}, volume No. 21 of {\em Van Nostrand Mathematical Studies}.
\newblock Van Nostrand Reinhold Co., New York-Toronto-London, 1969.

\bibitem{B1}
J.~J. Grobler.
\newblock Spectral theory in {B}anach lattices.
\newblock In {\em Operator theory in function spaces and {B}anach lattices}, volume~75 of {\em Oper. Theory Adv. Appl.}, pages 133--172. Birkh\"{a}user, Basel, 1995.

\bibitem{Example10}
A.~Guillin, B.~Nectoux, and L.~Wu.
\newblock Quasi-stationary distribution for {H}amiltonian dynamics with singular potentials.
\newblock {\em Probab. Theory Related Fields}, 185(3-4):921--959, 2023.

\bibitem{Yaglom2}
B.~Haas and V.~Rivero.
\newblock Quasi-stationary distributions and {Y}aglom limits of self-similar {M}arkov processes.
\newblock {\em Stochastic Process. Appl.}, 122(12):4054--4095, 2012.

\bibitem{kolb2}
G.~Hinrichs, M.~Kolb, and V.~Wachtel.
\newblock Persistence of one-dimensional {$\rm AR(1)$}-sequences.
\newblock {\em J. Theoret. Probab.}, 33(1):65--102, 2020.

\bibitem{Kolb1}
M.~Kolb and D.~Steinsaltz.
\newblock Quasilimiting behavior for one-dimensional diffusions with killing.
\newblock {\em Ann. Probab.}, 40(1):162--212, 2012.

\bibitem{UK}
U.~Krengel.
\newblock {\em Ergodic theorems}, volume~6 of {\em De Gruyter Studies in Mathematics}.
\newblock Walter de Gruyter \& Co., Berlin, 1985.

\bibitem{kr}
E.~Kreyszig.
\newblock {\em Introductory functional analysis with applications}.
\newblock Wiley Classics Library. John Wiley \& Sons, Inc., New York, 1989.

\bibitem{Yaglom1}
T.~G. Kurtz and S.~Wainger.
\newblock The nonexistence of the {Y}aglom limit for an age dependent subcritical branching process.
\newblock {\em Ann. Probability}, 1:857--861, 1973.

\bibitem{R1}
A.~Lambert.
\newblock Quasi-stationary distributions and the continuous-state branching process conditioned to be never extinct.
\newblock {\em Electron. J. Probab.}, 12:no. 14, 420--446, 2007.

\bibitem{CLM}
A.~Lasota and M.~C. Mackey.
\newblock {\em Chaos, fractals, and noise: stochastic aspects of dynamics}, volume~97.
\newblock Springer Science \& Business Media, 2013.

\bibitem{L}
J.~B. Lasserre and C.~E.~M. Pearce.
\newblock On the existence of a quasistationary measure for a {M}arkov chain.
\newblock {\em Ann. Probab.}, 29(1):437--446, 2001.

\bibitem{Lee2012IntroductionManifolds}
J.~M. Lee.
\newblock {\em Introduction to smooth manifolds}, volume 218 of {\em Graduate Texts in Mathematics}.
\newblock Springer, New York, second edition, 2013.

\bibitem{new}
H.~Li, H.~Zhang, and S.~Liao.
\newblock On quasi-stationaries for symmetric {M}arkov processes.
\newblock {\em J. Math. Anal. Appl.}, 528(1):Paper No. 127498, 18, 2023.

\bibitem{DifR5}
P.~Mandl.
\newblock Spectral theory of semi-groups connected with diffusion processes and its application.
\newblock {\em Czechoslovak Mathematical Journal}, 11(4):558--569, 1961.

\bibitem{Ideal1}
P.~Meyer-Nieberg.
\newblock {\em {B}anach lattices}.
\newblock Springer Science \& Business Media, 2012.

\bibitem{Tweedie}
S.~P. Meyn and R.~L. Tweedie.
\newblock {\em {M}arkov chains and stochastic stability}.
\newblock Springer Science \& Business Media, 2012.

\bibitem{milnor}
J.~Milnor.
\newblock {\em Dynamics in one complex variable}, volume 160 of {\em Annals of Mathematics Studies}.
\newblock Princeton University Press, Princeton, NJ, third edition, 2006.

\bibitem{Esperanca}
W.~O\c{c}afrain.
\newblock Quasi-stationarity and quasi-ergodicity for discrete-time {M}arkov chains with absorbing boundaries moving periodically.
\newblock {\em ALEA Lat. Am. J. Probab. Math. Stat.}, 15(1):429--451, 2018.

\bibitem{Yaglom3}
W.~O\c{c}afrain.
\newblock Polynomial rate of convergence to the {Y}aglom limit for {B}rownian motion with drift.
\newblock {\em Electron. Commun. Probab.}, 25:Paper No. 35, 12, 2020.

\bibitem{DifR6}
R.~G. Pinsky.
\newblock On the convergence of diffusion processes conditioned to remain in a bounded region for large time to limiting positive recurrent diffusion processes.
\newblock {\em Ann. Probab.}, 13(2):363--378, 1985.

\bibitem{survey}
P.~K. Pollett.
\newblock Quasi-stationary distributions: a bibliography.
\newblock {\em Available at \url{https://people.smp.uq.edu.au/PhilipPollett/papers/qsds/qsds.pdf}}, 2008.

\bibitem{RW}
L.~C.~G. Rogers and D.~Williams.
\newblock Diffusions, {M}arkov processes and martingales, volume 1: Foundations.
\newblock {\em John Wiley \& Sons, Ltd., Chichester}, 7, 1994.

\bibitem{rudin}
W.~Rudin.
\newblock {\em Real and complex analysis}.
\newblock McGraw-Hill Book Co., New York, third edition, 1987.

\bibitem{DifR8}
D.~Steinsaltz and S.~N. Evans.
\newblock Quasistationary distributions for one-dimensional diffusions with killing.
\newblock {\em Trans. Amer. Math. Soc.}, 359(3):1285--1324, 2007.

\bibitem{FiniteR4}
E.~A. van Doorn and P.~K. Pollett.
\newblock Quasi-stationary distributions for reducible absorbing {M}arkov chains in discrete time.
\newblock {\em Markov Process. Related Fields}, 15(2):191--204, 2009.

\bibitem{van}
Erik~A. van Doorn.
\newblock Quasi-stationary distributions and convergence to quasi-stationarity of birth-death processes.
\newblock {\em Adv. in Appl. Probab.}, 23(4):683--700, 1991.

\bibitem{V2}
D.~Villemonais.
\newblock {\em Exponential convergence to a quasi-stationary distribution and applications}.
\newblock Habilitation thesis, Universit{\'e} de Lorraine (Nancy), \textit{Available at} \url{https://www.normalesup.org/~villemonais/publication/hdr/}, 2019.

\bibitem{QEM1}
J.~Zhang, S.~Li, and R.~Song.
\newblock Quasi-stationarity and quasi-ergodicity of general {M}arkov processes.
\newblock {\em Sci. China Math.}, 57(10):2013--2024, 2014.

\bibitem{Ale}
H.~Zmarrou and A.~J. Homburg.
\newblock Bifurcations of stationary measures of random diffeomorphisms.
\newblock {\em Ergodic Theory Dynam. Systems}, 27(5):1651--1692, 2007.

\end{thebibliography}

\end{document}